\documentclass[11pt,a4paper]{amsart}
\usepackage{amsmath,caption,booktabs,lipsum}
\numberwithin{equation}{section}
\usepackage{amsthm}
\usepackage{pdfpages}
\usepackage{relsize}
\usepackage{amsfonts}
\usepackage{amssymb}
\usepackage{mathrsfs}
\usepackage{tikz}
\usepackage{caption}
\usepackage{environ}
\usepackage{elocalloc}
\usepackage{url}
\usepackage{caption}
\usepackage{graphicx}
\usepackage[backend = bibtex, style=numeric-comp,doi=false,isbn=false,url=false]{biblatex}
\usepackage{CJKutf8}
\usepackage[normalem]{ulem}
\usepackage{xcolor}
\usepackage{etoolbox}
\usepackage{caption}
\usepackage{mathtools}
\usepackage{dcpic}
\usepackage{tikz-cd}
\usepackage[bottom]{footmisc}
\usepackage{hyperref}
\usepackage{enumitem}
\usepackage{stmaryrd}
\usetikzlibrary{positioning}
\usetikzlibrary{arrows,chains,positioning,scopes,quotes}
\usetikzlibrary{decorations.markings}
\newcommand{\xn}{X\times N_0}
\newcommand{\yn}{Y\times N_f}

\newcommand{\supp}{\operatorname{Supp}}
\newcommand{\reg}{{\operatorname{Reg}}}
\newcommand{\ai}{\alpha}

\newcommand{\sing}{{\operatorname{Sing}}}

\newcommand{\mff}{\mathcal{F}}

\newcommand{\be}{\beta}

\newcommand{\ga}{\gamma}

\newcommand{\de}{\delta}

\newcommand{\e}{\epsilon}
\newcommand{\hm}{\mathcal{H}}
\newcommand{\lam}{\lambda}

\newcommand{\si}{\sigma}
\newcommand{\Si}{\Sigma}

\newcommand{\ta}{\theta}

\newcommand{\cms}{\operatorname{comass}}

\newcommand{\ms}{\mathbf{M}}

\newcommand{\cd}{\cdots}

\newcommand{\s}{\subset}

\newcommand{\no}[1]{\left\lVert#1\right\rVert}
\newcommand{\cu}[1]{\left\llbracket#1\right\rrbracket}

\newcommand{\dvol}{\operatorname{dvol}}

\newcommand{\du}{^\ast}
\newcommand{\pf}{_\ast}

\newcommand{\ka}{\kappa}
\newcommand{\pxn}{\pi_{X\times N_0}}
\newcommand{\pyn}{\pi_{Y\times N_f}}
\newcommand{\dvx}{\dvol_{X\times N_0}^{g_N}}
\newcommand{\dvy}{\dvol_{Y\times N_f}^{g_N}}
\newcommand{\m}{^{-1}}

\newcommand{\cz}{\C^{d-\dim N}/\Z^{2(d-\dim N)}}

\newcommand{\w}{\wedge}

\newcommand{\pd}{\partial}

\newcommand{\na}{\nabla}

\newcommand{\R}{\mathbb{R}}
\newcommand{\Z}{\mathbb{Z}}

\newcommand{\C}{\mathbb{C}}

\newcommand{\id}{\textnormal{id}}

\makeatletter
\def\thm@space@setup{%
	\thm@preskip=0.2cm plus 0cm minus 0cm
	\thm@postskip=\thm@preskip % or whatever, if you don't want them to be equal
}
\makeatother
\theoremstyle{plain}
\newtheorem{thm}{Theorem}
\newtheorem{exam}{Example}[section]
\newtheorem{lem}[exam]{Lemma}

\newtheorem{fact}[exam]{Fact}

\newtheorem{assump}[exam]{Assumption}

\newtheorem{defn}[exam]{Definition}

\theoremstyle{definition}
\newtheorem{conj}{Conjecture}

\title[Fractal singular sets]{On a conjecture of Almgren: area-minimizing submanifolds with fractal singular sets}
\author{Zhenhua Liu}
\dedicatory{Dedicated to Xunjing Wei}
\numberwithin{equation}{subsection}
\begin{document}
	\setlength{\abovedisplayskip}{5pt}
	\setlength{\belowdisplayskip}{5pt}
	\setlength{\abovedisplayshortskip}{5pt}
	\setlength{\belowdisplayshortskip}{5pt}
	\maketitle\vspace{-3em}
	\begin{abstract}
	We construct area-minimizing submanifolds with fractal singular sets on compact Riemannian manifolds. Thus, we settle a conjecture of Almgren and our answer is sharp dimension-wise. Furthermore, we can prescribe arbitrarily the strata in the Almgren stratification of the singular sets of our area-minimizing submanifolds, and our results hold in the category of integral currents, mod $v$ currents and stable stationary varifolds. 
	\end{abstract}
	\section{Introduction}
		The general problem of finding area-minimizing representatives in integral homology classes is solved by Federer and Fleming (\cite{FF}):
	\begin{quote}
		\emph{Every integral homology class on a compact Riemannian manifold admits an area-minimizing representative.}
	\end{quote}
	In other words we can always find a representative of an integral homology class that has the least area among all representatives of the same homology class. The representative in the above result is found in the category of integral currents, which roughly speaking is the closure of the set of algebraic topological polyhedron chains under Whitney flat topology \cite{FF}.
	
	Calling the area-minimizing representatives submanifolds was justified due to the following theorem:
	\begin{quote}
		\emph{An area-minimizing representative of an integral homology class is a smooth minimal submanifold counted with integer multiplicities outside of a singular set that is codimension $2$ countably rectifiable with respect to the submanifolds.}
	\end{quote}Almgren \cite{FA} first established that the singular set is of Hausdorff codimension $2,$ while the hypersurface case dated to  Federer \cite{HFts}. De Lellis and Spadaro simplified Almgren's proof in \cite{DS1,DS2,DS3}. Simultaneously and independently, De Lellis-Minter-Skoborotova 
	(\cite{DMS1,DMS2,DMS3}) and Wickramasekera-Krummel (\cite{BK1,BK2,BK3})  proved the above tour de force result. 
  
  Thus, we will use interchangeably the terms area-minimizing representatives of homology classes and area-minimizing submanifolds. Similar results hold for area-minimizing representatives of finite coefficient homology (\cite{DMSmod,LScy} and the references therein.)
	
	The above results on the singular sets of area-minimizing submanifolds are sharp, as singular sets do appear in very natural settings \cite[Section 4]{MR0168727}:
	\begin{quote}
		\emph{Complex algebraic subvarieties, including singular subvarieties, are area-minimizing in the Fubini-Study metric on complex projective spaces.}
	\end{quote}
In the same vein, many singular complex hyperquadrics are area-minimizing in every finite coefficient homology class they belong to (\cite{FMcalv}). 

To the author's knowledge, like complex algebraic subvarieties, all explicitly known examples of area-minimizing representatives of integral and finite coefficient homologies are subanalytic subvarieties. So are their singular sets. This is in sharp contrast with the above state of the art on the singular sets of area-minimizing currents, which gives us only the rectifiability of singular sets.

Already noticing this gap between general theories and known examples in the 1980s, Professor Frederick Almgren has raised the following question as Problem 5.4 in \cite{GMT}, (to quote him), 
\begin{conj}\label{ca}
		\emph{"Is it possible for the singular set of an area-minimizing integer (real) rectifiable current to be a Cantor type set with possibly non-integer Hausdorff dimension?"}
	\end{conj}
The ground-breaking work of Leon Simon \cite{LSfr} answered an analogue of the above question in the setting of stable stationary varifolds and inspired this manuscript (more comments in Section \ref{ls}). The main theorems of this manuscript answer the above conjecture by Almgren affirmatively and the answers are sharp in terms of dimensions of the singular sets.
\begin{thm}\label{thmd}
	For any integer $d\ge 2$ and any nonnegative real number $0\le\ai\le d-2,$ there exists a smooth compact $(d+3)$-dimensional Riemannian manifold $M^{d+3}$, and  a $d$-dimensional area-minimizing integral current $\Si$ on $M$ such that
\begin{enumerate}
\item The singular set of $\Si$ is of Hausdorff dimension $\ai$.
\item A $d$-dimensional smooth calibration form $\phi$ on $M$ calibrates $\Si.$
\end{enumerate}
\end{thm}
Theorem \ref{thmd} is a natural corollary of a much more general result below. To state the result, let us start with some basic assumptions and definitions that we will use throughout this manuscript.

For a $d$-dimensional stationary varifold $T$, which is a much weaker notion than being area-minimizing,  the Almgren stratification of $ T$ (\cite{BWst}) is an ascending chain of closed subsets of the support of  $T$
\begin{align*}
\mathcal{S}_0T\s\mathcal{{S}}_1T\s\cd\s\mathcal{S}_dT=\supp T,\end{align*}
such that $\mathcal{S}_j T$ consists of points in $\supp T$ with tangent cones having at most $j$-dimensional translational invariance. The reader can just regard the Almgren stratification as a geometric measure theory analogue of Whitney stratification.

Define the $j$-th stratum of the Almgren stratification of $T$ as points in $\supp T$ who has a tangent cone with precisely $j$-dimensional translation invariance. In other words,
\begin{defn}
We call $\mathcal{S}_jT\setminus\mathcal{S}_{j-1}T$ the $j$-th stratum of the Almgren stratification, with $S_{-1}T=\emptyset.$
\end{defn}
For integers $j\le d-1,$ points in the $j$-th stratum of $T$ necessarily lie in the singular set of $T.$ Now we are ready to state our assumptions, which roughly speaking are the data we will use in the Almgren stratification.
\begin{assump}\label{assumpmff}Assume the following
\begin{itemize}
\item	The symbol $d$ denotes an integer at least $2$ and the symbol $\mff$ denotes a finite collection of pairwise disjoint, smooth, closed, connected Riemannian manifolds, each of distinct dimension.
\item We set $k=\mathbf{min}_{N\in \mff}\dim N.$ For Theorem \ref{thmi} and Theorem \ref{thmp}, assume $\mathbf{max}_{N\in \mff}\dim N\le d-2$. For Theorem \ref{thms}, assume $\mathbf{max}_{N\in \mff}\dim N\le d-1.$
\item For each element $N\in \mff$, prescribe a closed subset $K_N\s N.$
\end{itemize}
\end{assump}
For instance, with $d=4,$ we can set $\mff=\{S^1,S^2\}$ and $K_{S^1},K_{S^2}$ Cantor sets in $S^1,S^2.$% We want to emphasis that elements of $\mathcal{F}$ will always appear as pairwise disjoint submanifolds in our constructions, so they do not have any intrinsic relations.

Under the above assumption, we can construct area-minimizing integral currents whose singular sets are the disjoint unions $\cup_{N\in \mff}K_N$ and each $K_N$ is precisely the $(\dim N)$-th stratum in the Almgren stratification. In other words, we can prescribe almost arbitrarily the structure of $T.$ 
\begin{thm}\label{thmi}Under Assumption \ref{assumpmff} there exists a smooth compact Riemannian manifold $M^{2d-k+1}$ of dimension $(2d-k+1)$, and  a $d$-dimensional area-minimizing integral current $\Si$ on $M$ such that
\begin{enumerate}
\item $\Si$ is the image of a smooth immersion into $M$. %$$M=\bigcup_{N\in\mff`} N\times (S^1)^{2d-k-\dim N+1},$$ 
\item The singular set of $\Si$ is the disjoint union $$\bigcup_{N\in\mff}K_N,$$   and the regular set of $\Si$ has finitely many connected components.
\item For all $N\in \mff $, the $(\dim N)$-th stratum of the Almgren stratification of the singular set of $\Si$ is precisely $K_N.$
\item A $d$-dimensional smooth calibration form $\phi$ on $M$ calibrates $\Si.$
\end{enumerate} 
\end{thm}
Recall that a set is countably $j$-rectifiable if it can be covered Hausdorff $j$-dimensional a.e. by countable unions of Lipschitz images of subsets of $\R^j.$ Thus,  Theorem \ref{thmi} above shows that the rectifiability of singular sets of area-minimizing integral currents established in \cite{DMS1,DMS2,DMS3,BK1,BK2,BK3} is sharp and cannot be improved much.

Both $\Si$ and $M$ in Theorem \ref{thmi} are not necessarily connected. Each stratum of the Almgren stratification of $\Si$ comes from the non-transverse intersection of two connected minimal submanifolds both lying inside one connected component of $M.$ Using techniques of Yongsheng Zhang (\cite{YZj,YZa}), we can make Theorem \ref{thmi} stronger by making both $\Si$ and $M$ connected. We will pursue this in future work \cite{ZLa2}. 

As a byproduct, we can also prove similar results for area-minimizing currents in $\Z/v\Z$-coefficient homology and stable stationary varifolds. 
\begin{thm}\label{thmp}
With the same notation as Thoerem \ref{thmi}, for every integer $v\ge 2$, the current $\Si$ is area-minimizing mod $v$ in the $\Z/v\Z$-coefficient homology class it belongs to.    Thus, any real number $0\le \ga\le d-2$ can be realized as the dimension of the singular set of a $d$-dimensional mod $v$ area-minimizing current.
\end{thm}
%Note that we do not claim we can make $\Si,M$ to have connected regular est in Theorem \ref{thmp}.
The above theorem in the case of $v=2$ is also sharp dimension-wise by \cite{HFts} and implies the sharpness of the rectifiability of singular sets proved in \cite{LScy}. For general $v$, the sharpness of the above theorem follows from \cite{DMSmod}.
\begin{thm}\label{thms}Under Assumption \ref{assumpmff} there exists a smooth compact Riemannian manifold $M^{2d-k+1}$ of dimension $(2d-k+1)$ and a stable stationary varifold $\Si$ in $M$ such that
\begin{itemize}
\item The varifold $\Si$ is the image of a smooth immersion.
\item The singular set of $\Si$ is the disjoint union $$\bigcup_{N\in\mff}K_N,$$ 
\item For all $N\in \mff $, the $(\dim N)$-th stratum of the Almgren stratification of the singular set of $\Si$ is precisely $K_N.$
\item For any diffeomorphism $\Phi$ that is smoothly homotopic to the identity, we have
\begin{align*}
	\ms(\Si)\le\ms(\Phi\pf\Si),
\end{align*}where $\ms$ denotes the mass of currents.
\end{itemize}  Thus, any real number $0\le \be\le d-1$ can be realized as the dimension of the singular set of a $d$-dimensional stable stationary varifold.
\end{thm}
The folklore conjecture is that $d$-dimensional stable stationary varifolds have a singular set of Hausdorff dimension at most $(d-1),$ so the above result is possibly sharp dimension-wise. The above theorem also establishes the sharpness of the main results in \cite{NV}.
\subsection{Plan of proof}\label{planpf}
Before stating the plan of our proof, we need a method to prove area-minimizing in general, and this method is via calibrations \cite{HL}.
	
	Recall that the comass of a  $d$-dimensional differential form  $\phi$ is the maximum of $\phi$ evaluated on unit simple $d$-vectors in the tangent space to $M$ among all  points \cite[Section 1.8]{HF}. In other words
	\begin{align}\label{eqcms}
		\cms_g\phi=\max_{q\in M}\max_{\substack{P\s T_qM\\ \dim P=d}}\phi(P).
	\end{align}
	Here we use the symbol $P$ to denote both an oriented $d$-dimensional plane $P$ and the unit simple $d$-vector representing $P.$
	\begin{defn}(Definition of calibrations)\label{defncal}
		\begin{itemize}
			\item 	We say a closed smooth $d$-form $\phi$ on a (possibly open) ambient manifold is a calibration if its comass is at most $1.$ 
			\item 
			We say a $d$-dimensional integral current $T$ is calibrated by $\phi,$ if $d$-dimensional Hausdorff measure almost everywhere $\phi$ restricted to the tangent space of $T$ equals the volume form of $T.$
		\end{itemize}
	\end{defn}
The fundamental theorem of calibrated geometry \cite[Theorem 4.2]{HL} gives	\begin{lem}\label{fcal}
		If a $d$-dimensional integral current $T$ is calibrated by a $d$-dimensional calibration form $\phi,$ then $T$ is area-minimizing.
		\end{lem}
The reader can try evaluating the complicated expression (\ref{eqcms}) with their favorite differential forms and will quickly realize that calculating the comass exactly is mission impossible. Thus, Frank Morgan has famously commented \cite[p. 343]{FMct}
\begin{quote}
\emph{Finding a calibration remains an art, not a science.}
\end{quote}

Since calibrations will be the main tool we use, the reader can imagine that our proof was obtained through trial and error with no guiding principles. Thus, the reader should keep in mind that our constructions will inevitably come out of the blue with not much heuristics.

We now outline the plan of our proof for the main theorems.
	
	Theorem \ref{thmd} is a direct corollary of Theorem \ref{thmi}, by taking $\mff=\{\R^{d-2}/\Z^{d-2}\}$ and the fact that a $(d-2)$-dimensional cube contains closed subsets of any Hausdorff dimension between $0$ and $(d-2).$ The latter classical fact follows by taking products of Cantor sets and intervals, e.g., \cite[Section 4.10 and Theorem 8.10]{PM}.

	From now on, we will focus on proving Theorems \ref{thmi}, \ref{thmp}, and \ref{thms}.
	
	Recall that our goal is to find an area-minimizing submanifold with fractal singular sets. Let us start with the simplest possible example of non-isolated singular sets and explain how we modify them.

Consider the standard $\C^3$ with a holomorphic coordinate system $(z_1,z_2,z_3)$ and the real coordinate system $(x_1,x_2,x_3,y_1,y_2,y_3)$ defined by $z_j=x_j+iy_j$ for $1\le j\le 3.$

By Lemma \ref{torf}, the form $dx_1dx_2dx_3+dy_1dy_2dx_3$ is a calibration form and calibrates the chain sum of the $x_1x_2x_3$-plane and the $y_1y_2x_3$-plane. Use the symbol $\cu{x_1x_2x_3}$ to denote the integral current corresponding to the $x_1x_2x_3$-plane and the symbol $\cu{y_1y_2x_3}$ to denote the integral current corresponding to the $y_1y_2x_3$-plane. Then the integral current
\begin{align*}
\cu{x_1x_2x_3}+\cu{y_1y_2x_3}
\end{align*} is area-minimizing by Lemma \ref{fcal}.%At this point, it is beneficial to give a robust explanation of why $dx_1dx_2dx_3+dy_1dy_2dx_3$ is a calibration form, as ultimately our construction depends on modifying $dx_1dx_2dx_3+dy_1dy_2dx_3$. For any $3$-plane $P$ in $\R^5,$ by diagonalization with respect to the orthogonal decomposition of $\R^5$ into $x_1x_2x_3$-plane and $y_1y_2$-plane \cite[Lemma 7.5]{HL}, $P$ is spanned by a set of orthonormal basis  $\{\cos \ta_1 v_1+\sin \ta_1 w_1,\,\cos\ta_2 v_2+\sin\ta_2 w_2,\,v_3\},$ with $v_1,v_2,v_3$ an orthonormal basis of the $x_1x_2x_3$-plane and $w_1,w_2$ an orthonormal basis of the $y_1,y_2$-plane, and $\ta_1,\ta_2\in(-\pi,\pi].$ Then we have \begin{align*}	&(dx_1dx_2dx_3+dy_1dy_2dx_3)(P)\\=&\cos\ta_1\cos\ta_2+\sin\ta_1\sin\ta_2dx_3(v_3)\\\le&|\cos\ta_1||\cos\ta_2|+|\sin\ta_1||\sin\ta_2|\\\le& |\cos(|\ta_1|-|\ta_2|)|.\end{align*}The key observation is that the two summands of $dx_1dx_2dx_3+dy_1dy_2dx_3$ are sufficiently orthogonal so that evaluating comass 

 The singular set of $\cu{x_1x_2x_3}+\cu{y_1y_2x_3}$ is the intersection set of $x_1x_2x_3$-plane with the $y_1y_2x_3$-plane, i.e., the $x_3$-axis.

Use $\Pi_{\cu{x_1x_2x_3}}$ to denote the nearest distance projection of $\R^6$ onto $\cu{x_1x_2x_3}$ and $\Pi_{\cu{y_1y_2x_3}}$ to denote the nearest distance projection of $\R^6$ onto $\cu{y_1y_2x_3}.$ Then we can interpret the form $dx_1dx_2dx_3+dy_1dy_2dx_3$ as the sum of pullbacks of volume forms of $\cu{x_1x_2x_3}$ and $\cu{y_1y_2x_3}$ under the respective projections, i.e.,
\begin{align*}
dx_1dx_2dx_3+dy_1dy_2dx_3=\Pi_{\cu{x_1x_2x_3}}\du \dvol_{\cu{x_1x_2x_3}}+\Pi_{\cu{y_1y_2x_3}}\du \dvol_{\cu{y_1y_2x_3}}\end{align*}

The idea to change the singular set of $\cu{x_1x_2x_3}
+\cu{y_1y_2x_3}$ is to push $\cu{y_1y_2x_3}$ parametrically to the $y_3$-direction. To be precise let $f:\R\to\R$ be a smooth function and define a parameterized integral current $\cu{y_1y_2f}$ by
\begin{align*}
\cu{y_1y_2f}=\{(0,0,x_3,y_1,y_2,f(x_3))|y_1,y_2,x_3\in\R\}.
\end{align*}

Then 
\begin{align*}
\Pi_{\cu{x_1x_2x_3}}\du \dvol_{\cu{x_1x_2x_3}}+\Pi_{\cu{y_1y_2f}}\du \dvol_{\cu{y_1y_2f}}
\end{align*}is a calibration form in a modified metric and calibrates $$\cu{x_1x_2x_3}
+\cu{y_1y_2f}.$$ Here $\Pi_{\cu{y_1y_2f}}$ is the nearest distance projection onto ${\cu{y_1y_2f}}$.
\begin{figure}
	\centerline{
		\includegraphics[width=0.9\paperwidth]{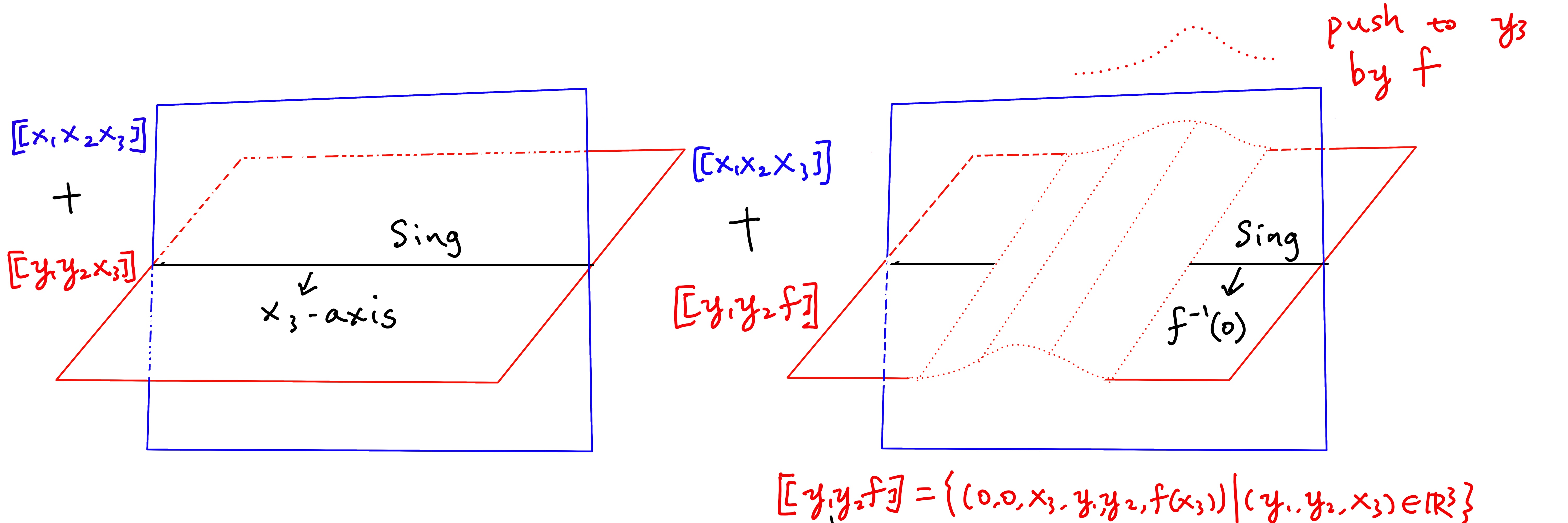}}
	\caption{Pushing $\cu{y_1y_2x_3}$ into $y_3$ direction parametrically by $f$}
	\label{figfractal}
\end{figure}
However, the singular set of the integral current
\begin{align*}\cu{x_1x_2x_3}
+\cu{y_1y_2f}
\end{align*}is the zero set of $f$ on $x_3$-axis. Thus, we can prescribe the singular set of $\cu{x_1x_2x_3}
+\cu{y_1y_2f}$ by prescribing $f$. The reader should refer to Figure \ref{figfractal} for intuition.

There is a hidden technicality in the above argument, as $\Pi_{\cu{y_1y_2f}}$ is not well defined everywhere. To make our construction work properly, we will extend $\Pi_{\cu{y_1y_2f}}$ to a retraction defined everywhere.

We can carry out the above construction in the torus $\C^3/\Z^6=\C^2/\Z^4\times S^1\times S^1$ and thus prescribe the first stratum in the Almgren stratification by requiring $f:S^1\to S^1$ vanishing to infinite order at its prescribed zero set in $S^1$. Other strata can be dealt with in a similar manner using the manifold
\begin{align}\label{mn}
\C^{d-\dim N}/\Z^{2(d-\dim N)}\times N\times S^1.
\end{align}
% To prove Theorem \ref{thmi}, we take Riemannian products of the manifolds in (\ref{mn}) with circles suitably many times to achieve the same ambient dimension and then take the disjoint union. To prove Theorem \ref{thmp}, we need to use calibrations modulo $v$ defined by \cite{FMcalv}. The proof of Theorem \ref{thms} relies on calibrating the connected components separately. The proof of Theorem \ref{thmic} uses Zhang's pioneering construction of calibrated connected sums \cite{YZa,YZj}. For the reader's convenience, we will fully reproduce Zhang's argument in this manuscript. 
	\subsection{Overview of the paper}
	Our paper will be structured as follows. 

	In Section \ref{bsdefn}, we will give some basic definitions and lemmas. 
	
In Section \ref{secss}, we will carry out our construction outlined in Section \ref{planpf} for a single stratum of the Almgren stratification.

In Section \ref{secpfmt}, we will prove that our main results follow from Section \ref{secss}.

%	In Appendices \ref{seccs} and \ref{seccsa}, we will reproduce the entire proof of the calibrated connect sum construction in Zhang's work \cite{YZa,YZj} and use them to prove Theorem \ref{thmic} in Section \ref{secpic}.
	
	In Section \ref{conc}, we will give some concluding remarks and discussions.
	\section*{Acknowledgements}
	The author acknowledges the support
of the NSF through the grant FRG-1854147.
I cannot thank my advisor Professor Camillo De Lellis enough for his unwavering support throughout the years. I would also like to thank him for giving me this problem and for many helpful discussions. Profound gratitude goes to Professor Leon Simon, whose groundbreaking work \cite{LS} has served as the primary inspiration for the author to pursue this problem. The author is deeply indebted to Professor Simon for his extraordinary generosity, guidance, and kindness throughout the years. On the occasion of his 80th birthday, the author warmly extends best wishes to Professor Simon.  For this manuscript, he has given many helpful suggestions, including the stronger form of Theorem \ref{thms} presented here. 

Sincere thanks go to Professor Hubert Bray, who has generously shared with the author invaluable advice about almost everything and whose words ultimately persuaded the author to pursue this conjecture of Almgren. The same thanks extend to Professors William Allard and Robert Bryant, who have taught me invaluable lessons on geometric measure theory and Riemannian geometry.  

The author is profoundly indebted to Professor Frank Morgan whose unwavering support and foundational works are always the main inspiration for the author. The heartfelt tributes extend to Professors Gary Lawlor, Dana Mckenzie and Yongsheng Zhang, whose works on calibrations form the core of this manuscript. Zhang's gluing construction of calibrations \cite{YZa,YZj} is what enables the author to extend local results to global ones in this work and many others.

Another thank you goes to Professors Marcus Khuri and Yang Li, whose comments about the calibrations used in a previous version of this manuscript inspired the author to find more general constructions used in this version. The same thank goes to Professors Kai Xu and Donghao Wang, who have helped the author navigate differential topology and geometry. The author would also like to thank Professors Simon Brendle, Antoine Song, and Song Sun for their constant encouragement and support.  Thanks go to the following colleagues for their interest in this work: Professors
	Nicolau Aiex, 
	Aidan Backus, Otis Chodosh, Fernando Cod\'{a} Marques, 
	Tobias Colding,
	 Tam\'{a}s Darvas, Benjamin Dees, Max Engelstein,  Herman Gluck, Or Hershkovits, Pei-Ken Hung, Robert Kusner, Chao Li, Yi Lai, 
Davi Maximo,	Peter McGrath,  Andr\'{e} Neves, Daniel Stern, and Ao Sun.

 Sincere thanks go to Professor David Kraines and my high school teacher Zhongyuan Dai, who have encouraged me to pursue mathematics. This thanks also extend to Professors Leslie Saper, Mark Stern and Shing-tung Yau, without whose help from the S.-T. Yau High School Science Award, I would not have been able to receive the wonderful mathematical education at Duke University as an undergraduate, where Professors Robert Calderbank,
  Richard Hain, and Colleen Robles have given constant encouragements.

Profound thanks go to the referees, whose careful reading of this manuscript and many generous and encouraging suggestions are what bring this manuscript to the current state. The author would also like to thank the referees for pointing out the reference \cite{MR1305283}.

My deepest thanks go to Dr. Lin Kong, Dr. Jiade Lu,  Dr. Paul Busse,  Dr. Michael Lanuti, Dr. Kartik Sehgal, and other members of my medical teams, without whose help I would not have survived cancer in 2019 and 2021. I apologize for the lengthy acknowledgment, but I cherish the opportunity to express my gratitude to those who have supported me along my journey. Facing significant morbidity risks in the coming years, I feel it is important to take this opportunity to formally acknowledge and thank all who have helped me.
	\section{Basic definitions}\label{bsdefn}
	We will give some basic definitions in this section. The contents are more or less standard, except for Sections \ref{tof}, \ref{secnb}, \ref{seccalv}, and the experienced reader can skip this section at the first reading.
\subsection{Euclidean spaces and flat tori}\label{seceuc}
For the standard $n$-dimensional complex Euclidean space $$\C^n,$$ we will always use a holomorphic coordinate system $$(z_1,\cd,z_n)$$ with the corresponding real coordinate system $$(x_1,\cd,x_n,y_1,\cd,y_n)$$ defined by 
$$z_1=x_1+iy_1,\cd,z_n=x_n+iy_n.$$ The standard flat metric on $\C^n$ is defined so that the collection of all real coordinate vector fields,  \begin{align}\label{bsxy}
\{\pd_{x_1},\cd,\pd_{x_n},\pd_{y_1},\cd,\pd_{y_n}\}
\end{align}
 form an orthonormal basis.

The symbol $$\C^n/\Z^{2n}$$ will be reserved for the tori obtained from the quotient of $\C^n$ by the lattice spanned by (\ref{bsxy}), i.e.,
\begin{align*}
\Z^{2n}=\textnormal{span}_{\Z}\{\pd_{x_1},\cd,\pd_{x_n},\pd_{y_1},\cd,\pd_{y_n}\}.
\end{align*}
The natural orthogonal direct sum splitting of 
\begin{align*}
	\C^n=\textnormal{span}_{\R}\{\pd_{x_1},\cd,\pd_{x_n}\}\oplus\textnormal{span}_{\R}\{\pd_{y_1},\cd,\pd_{y_n}\},
\end{align*}
gives a natural splitting of $\C^n$ as a Riemannian product
\begin{align}\label{cspl}
	\C^n=\cu{x_1\cd x_n}\times \cu{y_1\cd y_n},
\end{align}
where $\cu{x_1\cd x_n}$ denotes the $n$-dimensional oriented plane spanned by $\pd_{x_1},\cd,\pd_{x_n}$ and $\cu{y_1\cd y_n}$ denotes the $n$-dimensional oriented plane spanned by $\pd_{y_1},\cd,\pd_{y_n}.$

The product (\ref{cspl}) induces a natural Riemannian product splitting of
\begin{align*}
	\C^n/\Z^{2n}=X\times Y,
\end{align*}
where $X$ is the quotient of $\cu{x_1\cd x_n}$ by $\Z^{2n}$ and $Y$ is the quotient of $\cu{y_1\cd y_n}$ by $\Z^{2n}.$

The symbol $$S^1$$ will always mean the standard unit circle in $\R^2,$ regarded as a manifold. We will also frequently identify $S^1$ as the interval
\begin{align*}
[-{\pi}{},{\pi}{}],
\end{align*}
  with end points identified and we will reserve label $t$ to denote points in $S^1.$	\subsection{Manifolds}\label{bdmn}
	Though we deal with both integral and finite coefficient homologies in our manuscript, unless otherwise stated, every smooth manifold and smooth submanifold we mention in this manuscript is assumed to be orientable.  
	
	In general, we will use the symbol $d$ as the dimension of the integral currents and the symbol $c$ as the codimension of the current with respect to the ambient manifold. The symbol $M$ will be reserved for a $(d+c)$-dimensional closed ambient manifold and the symbol $\Si$ will be reserved for a $d$-dimensional integral current of codimension $c$ inside $M.$ %When we do not mention the ambient manifold explicitly, it is understood that the ambient manifold is $M.$ 
	
	We will often speak of smooth neighborhoods or smooth open sets, by which we mean an open set with smooth boundary.%	We will also make use of transversality \cite{MH}. By this we mean that whenever we have a  pair of smooth submanifolds of possibly different dimension, we can assume without loss of generality that the pair is transverse to each other and thus intersects along a submanifold.
	
		We will also use transversality \cite{MH} frequently. By this we mean that whenever we have a  pair of smooth submanifolds of dimensions $d,d'$ on an ambient manifold of dimension $(d+c),$ respectively, we can arrange by an arbitrarily small smooth perturbation so that the pair is transverse to each other and thus intersects along a submanifold of dimension $(d'-c)$.	
	\subsection{Riemannian geometry}\label{basrie}
We will use $\textnormal{dist}(p,K)$ to denote the Riemannian distance between a point $p$ and a closed set $K.$ 
	
	We will reserve the symbol $\no{\cdot}$ to denote the Riemannian length of vectors, forms, etc. Whenever we use the symbol $|\cdot|,$ it is understood that we are either taking absolute value or taking the square root of the sum of squares of components in coordinates, thus different from $\no{\cdot}$ on general manifolds.
\subsection{Planes, vectors and forms}
For a $d$-dimensional oriented plane $P$ at a point in the tangent space of the ambient manifold, we will also use the symbol $P$ to denote the unique simple unit $d$-vector generating $P$ and will use $P^\ast$ to denote the $d$-form Riemannian dual to $P.$ 

Conversely, when we write a unit simple $d$-vector $v_1\w\cd\w v_d$, it can also mean the $d$-dimensional oriented plane spanned by $v_1,\cd,v_d.$
\subsection{Several facts about comass}
With the definition of comass in Section \ref{planpf} in mind, we record some facts about comass on \textbf{vector spaces}. We will frequently apply the below fact \textbf{pointwise} in tangent spaces of ambient manifolds when constructing calibrations and calculating comass.
\begin{fact}\label{cmsvec}
	Let $g,h$ be positive definite quadratic forms on a finite-dimensional vector space, with $\psi$ a constant differential form, i.e., a form with constant coefficients everywhere with respect to a fixed dual basis.\begin{enumerate}
		\item $\cms_h\psi\le\cms_g\psi$ if $h\ge g$ as quadratic forms.\label{cms1}
		\item $\cms_{\lam^2 h}\psi=|\lam|^{-\dim\psi}\cms_h\psi$.\label{cms2}
			\item If $\phi$ is a simple form, then $\cms_g\phi=\no{\phi}_g$, i.e., its comass equal to its Riemannian length.\label{cms0}
	\end{enumerate}
\end{fact}
Here the parameter $\lam$ satisfies $\lam\in \R\setminus\{0\}$, the symbol
$\cms_g\psi$ means the comass of $\psi$ with respect to $g$. The first two bullets are proved in \cite[Lemma 2.1.9, 2.1.20]{YZt}. The last bullet is a direct calculation using Cauchy-Binet.\subsection{A torus differential form}\label{tof}
In this section, we collect a basic lemma about calibration form that will be central to our constructions.
\begin{lem}
	\label{torf}
	Assume that
	\begin{enumerate}
		\item $V,W,Z$ are three finite-dimensional vector spaces equipped with positive definite inner products $g_V,g_W,g_Z$ and satisfy $\dim_V=\dim_W\ge 2.$
		\item $\phi_V,\phi_W,\phi_Z,\phi_Z'$ are constant coefficient non-vanishing simple differential forms on $V,W,Z,Z$, respectively, such that 
		\begin{align*}
			\dim\phi_V=&\dim\phi_W\ge 2,\\\dim\phi_Z=&\dim\phi_Z'.
		\end{align*}
	\end{enumerate} 
	Then on the standard Riemannian product $$V\times W\times Z$$ equipped with the product metric $$g=g_V\circ\pi_V+g_W\circ\pi_W+g_Z\circ\pi_Z,$$ we have
	\begin{align}
	&	\cms_g \bigg(\pi_V\du\phi_V\w \pi_Z\du\phi_Z+\pi_W\du\phi_W\w\pi_Z\du\phi_Z'\bigg)\\=&\max\bigg\{\cms_{g}\big(\pi_V\du\phi_V\w \pi_Z\du\phi_Z\big),\cms_{g}\big(\pi_V\du\phi_V\w \pi_Z\du\phi_Z'\big)\bigg\}.\label{cmsc}
	\end{align}
\end{lem}
The above lemma is false for $\dim\phi_V=\dim\phi_W=1.$
\begin{proof}
First of all, by choosing orthonormal frames $\{v_1,\cd,v_{\dim V}\},\{w_1,\cd,w_{\dim W}\}$ on $V,W,$ respectively, we can assume that
\begin{align*}
	\phi_V=&\ai_Vv_1\du\w\cd\w v_{\dim\phi_V}\du,\\
	\phi_W=&\ai_W w_1\du\w\cd\w w_{\dim \phi_W}\du,
\end{align*}
for some $\ai_V,\ai_W\not=0.$

We use \cite[Lemma 2.1  with $e=v_1,f=v_2,\psi=\ai_Vv_3\w\cd\w v_{\dim \phi_V}\w\pi_Z\du \phi_Z,\chi=\pi_W\du\phi_W\w\pi_Z\du\phi_Z',$]{DHM}, we deduce
\begin{align*}
	& \cms_g \bigg(\pi_V\du\phi_V\w \pi_Z\du\phi_Z+\pi_W\du\phi_W\w\pi_Z\du\phi_Z'\bigg)\\=&\max\bigg\{\cms_{g}\big(\ai_Vv_3\du\w\cd\w v\du_{\dim \phi_V}\w\pi_Z\du \phi_Z\big),\cms_{g}\big(\pi_V\du\phi_V\w \pi_Z\du\phi_Z'\big)\bigg\}\\
=&\max\bigg\{\cms_{g}\big(\pi_V\du\phi_V\w \pi_Z\du\phi_Z\big),\cms_{g}\big(\pi_V\du\phi_V\w \pi_Z\du\phi_Z'\big)\bigg\},
	\end{align*}where at the last line we have used \cite[Proposition 7.10]{HL}\end{proof}
\subsection{Product of simple forms}
\begin{lem}\label{lemvfm}
		Assume that
	\begin{enumerate}
		\item $\mathcal{V}$ is a finite collection of finite-dimensional vector spaces equipped with positive-definite quadratic forms.
		\item For each element $(V,g_V)$ of $\mathcal{V},$ we associate a non-vanishing constant simple form $\phi_V$ on $V$, a non-vanishing constant simple multi-vector $\xi_V$, with $\dim\phi_V=\dim\xi_V$, and a positive real number $\ai_V>0.$
		\item We equip the direct product 	\begin{align*}
			X_{\mathcal{V}}=	\prod_{V\in\mathcal{V}}V=\bigoplus_{V\in\mathcal{V}}V,
		\end{align*}with the Riemannian metric
		\begin{align*}
		g_\ai=\sum_{V\in\mathcal{V}}\ai_V^{\frac{2}{\dim\phi_V}} g_V\circ\pi_V.
		\end{align*}
	\end{enumerate}
Then	for any non-empty ordered subset $\mathcal{V}'$ of $\mathcal{V},$ the constant form
\begin{align}\label{bwv}
	\bigwedge_{V\in \mathcal{V}'}\pi_V\du\phi_V,
\end{align}satisfies
\begin{itemize}
	\item The multi-vector normalization equality
	\begin{align}\label{eqn}
		\bigg(\bigwedge_{V\in \mathcal{V}'}\pi_V\du\phi_V\bigg)\bigg(\bigwedge_{V\in \mathcal{V}'}\frac{\xi_V}{\no{\xi_V}_{g_\ai}}\bigg)=\prod_{V\in \mathcal{V}'}\ai_V\m\frac{\phi_V(\xi_V)}{\no{\xi_V}_{g_V}}.
	\end{align}
	\item The comass equality
	\begin{align}\label{eqw}
		\cms_{g_\ai}\bigwedge_{V\in \mathcal{V}'}\pi_V\du\phi_V=\prod_{V\in\mathcal{V}'}\ai_V\m\cms_{g_V}\phi_V.
	\end{align}
\end{itemize}
\end{lem}
Here $\pi_V$ is the canonical projection of the product $X_{\mathcal{V}}$ onto its factor $V.$ Ordered subset $\mathcal{V}'$ means that we assign numbers $1,2,\cd,\textnormal{card}(\mathcal{V'})$ bijectively to elements of $\mathcal{V}'$, and the wedge products (\ref{bwv}) and (\ref{eqw}) are taken in this order.
\begin{proof}
Definition of wedge product gives$\bigg(\bigwedge_{V\in \mathcal{V}'}\pi_V\du\phi_V\bigg)\bigg(\bigwedge_{V\in \mathcal{V}'}{\xi_V}\bigg)=\prod_{V\in \mathcal{V}'}\phi_V(\xi_V).$ The equation (\ref{eqn}) follows.
	
By bullet (\ref{cms2}) of Fact \ref{cmsvec}, to prove (\ref{eqw}), it suffices to prove that
\begin{align}\label{eqcmsa}
	\cms_{g_\ai}	\bigwedge_{V\in \mathcal{V}'}\pi_V\du\phi_V=\prod_{V\in \mathcal{V}'}\cms_{g_\ai|_{V}}\phi_V,
\end{align}where $g_\ai|_{V}$ means the restriction of $g_\ai$ onto the subspace $V$ in the direct sum $X_{\mathcal{V}}$. The above (\ref{eqcmsa}) follows directly from the (\ref{cms0}) in Fact \ref{cmsvec} and the definition of Riemannian products.
\end{proof}
\subsection{The normal bundle calibration}\label{secnb}
Let $W$ be a closed $w$-dimensional smooth submanifold of a closed Riemannian manifold $(V,h).$

Suppose a compact set $U(W)$ containing $W$ lies in the bijective image of the normal bundle exponential map $\exp^\perp_W$ of $W.$
\begin{lem}\label{lemznb}
The closed $w$-form $$\Pi_W\du\dvol_W^h$$ calibrates $W$ in $U(W)$ with the smooth metric $$h'=\no{\Pi_W\du\dvol_W^h}_h^{\frac{2}{w}}h.$$
\end{lem}
Here $\Pi_W$ is the nearest distance projection onto $W,$ $\dvol_W^h$ is the volume form of $W$ in $h$ and $\no{\cdot}_h$ is the Riemannian length in $h.$
\begin{proof}
	Since $U(W)$ lies in the bijective image of the normal bundle exponential map, it is straightforward to verify that $\no{\Pi\du\dvol_W^h}_h$ has a positive lower bound on $U(W).$ Thus, $h'$ is smooth. Note that volume forms are always simple and closed. Since pullback of forms preserves simpleness and closedness, $\Pi\du\dvol_W^h$ is also simple and closed. By bullet (\ref{cms0}) of Fact \ref{cmsvec}, $\no{\Pi\du_W\dvol_W^h}$ equals the comass of $\Pi\du_W\dvol_W^h$ in tangent spaces to $V.$ Then our lemma is just reformulating \cite[Remark 3.5 and Lemma 3.4 ]{YZa}.	The reader can also verify the above lemma using direct calculations via bullets (\ref{cms1}) and (\ref{cms0}) of Fact \ref{cmsvec}.\end{proof}
\subsection{Definition of area-minimizing}
	The right category to discuss area-minimizing representatives is the category of integral currents and integral currents modulo $v$, with integer $v$ at least $2$.  We will use Federer's definitive monograph \cite{HF} and Simon's classical lecture notes \cite{LS} as basic references. 
	
	For our purposes, the reader can heuristically regard integral currents as integer coefficient simplicial chains in algebraic topology. In reality, integral currents can be defined \cite[4.1.24]{HF} as the closure of the set of simplicial chains under mass topology applied both to the chains and their boundaries. Similarly,  the reader can just regard integral currents mod $v$ as limits of $\Z/v\Z$-coefficient simplicial chains in algebraic topology. Note that integral currents and integral currents mod $v$ are allowed to have finite area boundaries.
	
	We say an integral current $T$ is area-minimizing, if $T$ has the least area among all integral currents homologous to $T$:
	\begin{defn}\label{defnam}
		An $d$-dimensional integral current $T$ is area-minimizing if 
		\begin{align*}
			\ms(T)\le \ms(T+\pd V),
		\end{align*}for all  $(d+1)$-dimensional integral currents $V$. 
	\end{defn}
	Here $\ms$ is the mass of the current, which in our context, is just the area of the underlying set with multiplicity included. For instance, $\ms(\pm 2 S^1)=2\ms(S^1)=4\pi$ for the unit circle $S^1$ in $\R^2.$

The primary tool to prove area-minimizing is calibrations, which the reader can recall from Section \ref{planpf} and the main reference is \cite{HL}.

Similarly, we say an integral current $T$ is area-minimizing mod $v$, if $T$ has the least area among all integral currents homologous to $T$ in $\Z/v\Z$ homology:
	\begin{defn}\label{defnamv}
		An $d$-dimensional integral current $T$ is area-minimizing mod $v$ if 
		\begin{align*}
			\ms(T)\le \ms(T+\pd W+vZ),
		\end{align*}for all  $(d+1)$-dimensional integral currents $W$ and $d$-dimensional integral currents $Z$. 
	\end{defn}
The above definition of mod $v$ area-minimizing for integral currents is similar to the one in \cite[Definition 1.2]{DPHM}. Strictly speaking, the canonical way to define mod $v$ area-minimizing is to invoke the notions of rectifiable currents, mass mod $v$ and flat chains mod $v$. However, by  \cite[Corollaries 1.5, 1.6]{RY} and denseness of the set of integral currents inside the set of rectifiable currents under mass norm, the above definition is equivalent to the canonical one in our uses. 
	\subsection{Definition of singular sets}
	\begin{defn}\label{defnsm}
		We say an integral current $T$ is smooth at a point $p$ in the support of $T$ if there exists an open set $U$ containing $p$ on $M$, such that  $T$ restricted to $U$ equals an integer multiple of a smooth submanifold $N.$ The definitions of regular sets and singular sets of $T$ are as follows:
		\begin{itemize}
			\item The singular set of $T,$ $\sing T,$ is defined as set of points in the support of $T$ where $T$ is not smooth.		
			\item 
			The regular set of $T,$ $\operatorname{Reg}T$, is defined as the set of points in the support of $T$ where $T$ is smooth. 
		\end{itemize}
	\end{defn} 
	Here support means the underlying set of an integral current. From now on, we will also use the symbol $\supp T$ to mean the support of $T.$ %, and unfortunately is the standard term in geometric measure theory. The reader should not confuse it with the notion of supports in simplicial complexes.
	%In other words a current $T$ is smooth at a point $p$ in the underlying set of $T$, if the current $T$ restricted to a neighborhood of $p$ equals $kN$ for an integer $k\in N$ and a smooth submanifold $N$. The singular set of $T$ is defined as the complement of the smooth points of $T$.In our manuscript, all integral currents we construct will be chains of simplicial complexes with explicit descriptions. It will always be clear from the explicit descriptions alone what are the singular set and the regular set of the currents we  construct, so we will not spend time on formally proving what subsets are the singular sets.
	For example, the figure $8$ has a singular point at its self-intersection, while figure $0$ counted with multiplicity $-2$ is smooth.
\subsection{Area-minimizing integral currents in Riemannian products and disjoint unions}\label{basdis}
With the definition of area-minimizing in the previous subsection, we will collect several facts about area-minimizing integral currents.
\begin{lem}\label{lemdis}
Suppose the $d$-dimensional integral currents $T_1,\cd,T_n$ are area-minimizing (in integral or mod $v$ homology) inside compact Riemannian manifolds $M_1,\cd,M_n$, respectively. Then the sum of currents
\begin{align*}
T_1+\cd+T_n
\end{align*}
is area-minimizing (integral or mod $v$, respectively) in the disjoint union manifold
\begin{align*}
M_1\cup\cd \cup M_n.
\end{align*}
Furthermore, if $T_1,\cd,T_n$ are calibrated by $\phi_1,\cd,\phi_n$ on $M_1,\cd,M_n$, respectively, then $T_1+\cd+T_n$ is calibrated by the unique closed form whose restriction to $M_1,\cd,M_n$ is $\phi_1,\cd,\phi_n$, respectively. 
\end{lem}
\begin{proof}
Let $T$ be an area-minimizing integral current in $M_1\cup \cd\cup M_n$ homologous to $T_1+\cd+T_n$. Since $T$ is area-minimizing, we deduce that 
\begin{align}\label{mstn1}
\ms(T)\le\ms(T_1)+\cd+\ms(T_n).
\end{align}
On the other hand, the homology of $M_1\cup \cd\cup M_n$ canonically decompose as a direct sum
\begin{align}\label{homdec}
H_\ast(M_1\cup \cd\cup M_n)=H_\ast(M_1)\oplus\cd\oplus H_\ast(M_n).
\end{align}
Thus, $T$ restricted to $M_1,\cd,M_n$, must be homologous to $T_1,\cd,T_n$ inside $M_1,\cd,M_n$, respectively. From the area-minimality of $T_1,\cd,T_n$ inside $M_1,\cd,M_n$, respectively we deduce that  
\begin{align}\label{mstn2}
\ms(T)\ge\ms(T_1)+\cd+\ms(T_n).
\end{align}
Combining (\ref{mstn1}) and (\ref{mstn2}), we deduce that $\ms(T)=\ms(T_1)+\cd+\ms(T_n)$ and thus $T_1+\cd+T_n$ is also area-minimizing.

The claim about calibration follows directly from definitions.
\end{proof}
\begin{lem}\label{lemprod}
Suppose the $d$-dimensional integral current $T$ is area-minimizing (in integral or mod $v$ homology) on a compact  Riemannian manifold $V$. Let $W$ be another  compact Riemannian manifold with $p$ a point in $W$. Define a map $$i_p:V\to V\times W,$$ by
\begin{align*}
i_p(q)=(q,p).
\end{align*}
Then the pushforward current 
\begin{align*}
(i_p)\pf(T)
\end{align*}
is area-minimizing (in integral or mod $v$ homology) on $V\times W$. Furthermore, if $T$ is calibrated by a smooth form $$\phi$$ on $V,$ then $(i_p)\pf T$ is calibrated by a smooth form
\begin{align*}
\pi_V\du\phi,
\end{align*}on $V\times W,$ where $\pi_V$ is the canonical projection of $V\times W$ onto the $V$ factor.
\end{lem}
\begin{proof}
Let us first note that both $\pi_V$ and $i_q\circ\pi_V$ are area-non-increasing. In other words, if $\nu$ is a simple $d$-vector on the tangent space to $V\times W,$ then we have 
\begin{align}
&\no{(\pi_V)\pf \nu}\le \no{\nu},\label{pman}.% \\&\no{(i_q\circ\pi_V)\pf \nu}\le \no{\nu}\label{pqan}. 
\end{align}
Inequality (\ref{pman}) follows directly from the definition of Riemannian length of $d$-vectors and Cauchy-Binet formula. %Inequality (\ref{pqan}) follows from (\ref{pman}) and the fact that $i_p$ is a Riemannian isometric embedding.

Now let $Z$ be an integral current homologous to $(i_p)\pf T$ in $V\times W.$ Then $(\pi_V)\pf Z$ is homologous to $( \pi_V\circ i_p)\pf T$ on $V$.  By definition, $ \pi_V\circ i_p$ is the identity map on $V.$ Thus, $(\pi_V)\pf Z$ is homologous to $T.$ By the area-minimality of $T,$ we have
\begin{align*}
\ms(T)\le \ms((\pi_V)\pf Z).
\end{align*}
By (\ref{pman}), we have
\begin{align*}
\ms((\pi_V)\pf Z)\le\ms(Z).
\end{align*}
Finally, since $i_p$ is a Riemannian isometry, we deduce that
\begin{align*}
\ms((i_p)\pf T)=\ms(T)\le \ms((\pi_V)\pf Z)\le\ms(Z).
\end{align*}
Thus, $(i_p)\pf T$ is area-minimizing.

Direct calculation using (\ref{eqcms}) and (\ref{pman}) shows that $\pi_V\du \phi$ has comass at most that of $\phi$ and $(\pi_V)\du\phi$ calibrates $(i_p)\pf T.$ 
\end{proof}
\subsection{Calibration mod $v$ and mod $v$ area-minimizing}\label{seccalv}
In this subsection, we will prove a lemma about using special calibrations to prove mod $v$ area-minimizing. The following lemma is essentially due to Frank Morgan \cite{FMcalv}.
\begin{assump}\label{assumpv}Assume that
	\begin{enumerate}		
		\item  $M$ is a $(d+c)$-dimensional Riemannian manifold
		\item  $\mathcal{W}$ is a finite collection of closed connected $d$-dimensional oriented submanifolds of $M$.
		\item  For each element $W\in\mathcal{W},$ we have a smooth retraction: $\pi_W:M\to W$ onto $W.$
		\item \label{sigsum} For any map $\si:\mathcal{W}\to\{-1,1\},$ the form
		\begin{align*}
			\sum_{W\in\mathcal{W}}\si(W)\pi_W\du\dvol_W
		\end{align*}is a calibration form on $M$.
	\end{enumerate}
\end{assump}
\begin{lem}\label{lemcalv}Under Assumption \ref{assumpv}, if
$\Si$ is a $d$-dimensional integral current cycle on $M$, calibrated by the form, $$\sum_{W\in\mathcal{W}}\pi\du_W\dvol_W,$$
and for every $w\in \mathcal{W},$ we have
\begin{align*}
	(\pi_W)\pf \Si=u_WW
\end{align*}for some integer $$-\frac{v}{2}\le u_W\le \frac{v}{2},$$
then $\Si$ is area-minimizing mod $v$ (Definition \ref{defnamv}).
\end{lem}
For instance any degree $d$ holomorphic subvariety of $\C^n/\Z^{2n}$ is area-minimizing mod $v$ for $v\ge 2d$ (\cite{FMcalv}).
\begin{proof}
 Recall Definition \ref{defnamv}. Let $Z$ be any $d$-dimensional integral current on $M$ and $L$ be any $(d+1)$-dimensional integral current on $M.$ We need to prove that
\begin{align}
\label{ineq0}	\ms(\Si)\le\ms(\Si+\pd L+vZ).
\end{align}
We will break this down into three successive inequalities.
	\begin{align}
\label{ineq1}		\ms(\Si+\pd L+vZ)\ge&\sum_{W\in\mathcal{W}}\ms\bigg((\pi_W)\pf \big(\Si+\pd L+vZ\big)\bigg),\\
		\label{ineq2}\ms\bigg((\pi_W)\pf \big(\Si+\pd L+vZ\big)\bigg)\ge &\ms((\pi_W)\pf \Si),\textnormal{ for all }W\in\mathcal{W}\\
		\label{ineq3}\sum_{W\in\mathcal{W}}\ms((\pi_W)\pf \Si)=&\ms(\Si).			\end{align}
Before starting the proof, it is beneficial to review the integration representatives of integral currents. For a $d$-dimensional integral current $T$, integrating a $d$-dimensional differential form $\phi$ on $T$ or on $\Pi\pf T$, with $\Pi$ a smooth map, or calculating the mass of $T$ can all be written as integrals   (\cite[Theorem 4.1.28 (4)]{HF})
\begin{align}
	T(\phi)=&\int_{\supp T}\ta_T \phi(\nu_T)d\hm^d,\\
	\Pi\pf	T(\phi)=&\int_{\supp T}\ta_T \Pi\du\phi(\nu_T)d\hm^d,\\
	\label{eqint}	\ms(T)=&\int_{\supp T}\ta_T d\hm^d.
\end{align}
Here $\supp T$ is the support of $T$, which the reader can just assume to be the underlying set of $T$ and $d\hm^d$ means integrating with respect to $d$-dimensional Hausdorff measure. The symbol $\nu_T$ means the oriented tangent planes of $T$, normalized to be a unit simple $d$-vector, which exist Hausdorff $d$-dimensional everywhere and the symbol $\ta_T$ means the density of $T$, which is just the positive integer multiple we assign to $\nu_T$. Alternatively, the mass of $T$, i.e., the area of $T$, can be defined as
(\cite[p. 358 Section 4.1.7]{HF})
\begin{align}\label{eqintm}
	\ms(T)=\sup_{\textnormal{smooth }d\textnormal{-form }\phi\textnormal{ with }\cms\phi\le 1}T(\phi).
\end{align}
If the ambient manifold $T$ resides in is $d$-dimensional and oriented, then straightforward calculation shows that
\begin{align}	
\label{eqdst}	\nu_T=&\pm \nu_M,\\
\label{eqintma}	\ms(T)=&\sup_{f\in C^\infty(M,[-1,1])}T(f\dvol_M),
\end{align}
where $\nu_M$ is the oriented tangent space to $M.$
		
Now we are ready to do the proof, for (\ref{ineq1}), using the integration representative (\ref{eqint}) of a $d$-dimensional integral $T$ on $M$ and (\ref{sigsum}) in Assumption \ref{assumpv}, we have
\begin{align}
	&\ms(T)\\=&\int_{\supp(T)}\ta_{T}d\mathcal{H}^d\\
	\ge&\int_{\supp(T)}\ta_{T}\textnormal{ }\bigg(\sum_{W\in\mathcal{W}}|\pi\du_W\dvol_W(\nu_{T})|\bigg)d\mathcal{H}^d\\
	\ge&\sum_{W\in\mathcal{W}}\sup_{f\in C^\infty(W,[-1,1])}\int_{\supp(T)}\ta_{T}\textnormal{ }\pi\du_W(f\dvol_W)(\nu_{T})d\mathcal{H}^d\\
	=&\sum_{W\in\mathcal{W}}\ms((\pi_W)\pf (T)),\label{stepf}
\end{align}
where at (\ref{stepf}) we have used (\ref{eqintma}).

Applying $T=\Si+\pd L+vZ$ and $T=\Si$ yields
\begin{align}
	\ms(\Si+\pd L+vZ)\ge& \sum_{W\in\mathcal{W}}\ms((\pi_W)\pf(\Si+\pd L+vZ)),\\
\label{ineqms}\ms(\Si)\ge&	\sum_{W\in\mathcal{W}}\ms((\pi_W)\pf (\Si)).
\end{align}
In other words, we have proved (\ref{ineq1}).

For (\ref{ineq2}), we have
\begin{align}\label{eqpw}
	(\pi_W)\pf(\Si+\pd L+vZ)=u_WW+\pd((\pi_W)\pf L)+v(\pi_W)\pf Z=u_WW+v(\pi_W)\pf Z,
\end{align}
where we have used the fact that there are no non-trivial $(d+1)$-dimensional currents on a $d$-dimensional manifold, simply by definition of currents.

Since $(\pi_W)\pf(\Si+\pd L+vZ)$ is a $d$-dimensional integral current on the $d$-dimensional closed manifold $M,$ writing in integration form, by (\ref{eqpw}), we have\begin{align}
	&\ms((\pi_W)\pf(\Si+\pd L+vZ))\\=&\int_{\supp(u_WW+v(\pi_W)\pf Z)}\ta_{u_WW+v(\pi_{W})\pf Z}d\mathcal{H}^d\\
\label{stp0}	=&\int_{\supp W\setminus\supp (\pi_W)\pf Z}|u_W|d\mathcal{H}^d+\int_{\supp(\pi_W)\pf Z}|u_W\pm v\ta_{(\pi_W)\pf Z}|d\mathcal{H}^d\\
\label{stp0.5}	\ge&\int_{\supp W\setminus\supp (\pi_W)\pf Z}|u_W|d\mathcal{H}^d+\int_{\supp(\pi_W)\pf Z}\frac{v}{2}d\mathcal{H}^d\\
\label{stp1}	\ge&\int_{W}|u_W|d\mathcal{H}^d\\
	=&\ms((\pi_W)\pf \Si).
\end{align}
At (\ref{stp0}), the $\pm$ sign means either $+$ or $-$, with possibly different choices at each point of integration. The sign comes from (\ref{eqdst}). At (\ref{stp0.5}) and (\ref{stp1}), we have used $|u_W|\le\frac{v}{2}$ and the fact that density of $d$-dimensional integral currents are at least $1$ $d$-dimensional Hausdorff everywhere in its support. Thus $|u_w\pm v\ta_{(\pi_W)\pf Z}|\ge \frac{v}{2}\ge |u_W|.$
%By the constancy theorem \cite[Theorem 4.9]{FMgmt} and $\pd \Si=0,$ we always have $(\pd_W)\pf\Si=u_WW$ for some $u_W\in\Z.$ By our assumption, we always have $|u_W|\le \frac{v}{2}.$

For (\ref{ineq3}), %by \cite[Lemma 2.3 and proof of main theorem]{FMcalv}, at $d$-dimensional Hausdorff measure almost every point in the support of $\Si,$ for all $w\in\mathcal{W},$ we always have \begin{align}\label{eqfm}	\pi_W\du\dvol_W(\nu_\Si)\ge0.\end{align} This implies that
by (\ref{eqintma}) we have
\begin{align}\label{eqml}
\Si(\pi_W\du\dvol_W)=(\pi_W)\pf\Si(\dvol_W)\le\ms((\pi_W)\pf \Si).
\end{align}  Since $\Si$ is calibrated by $\sum_{w\in\mathcal{W}}\pi_W\du\dvol_W,$ i.e., $$\sum_{w\in\mathcal{W}}\pi_W\du\dvol_W(\nu_\Si)=1,$$ by (\ref{eqml})  we deduce that
\begin{align*}
	\ms(\Si)=	\sum_{w\in\mathcal{W}}\Si(\pi\du_W\dvol_W)\le\sum_{w\in\mathcal{W}}\ms((\pi_W)\pf \Si).
\end{align*}
Using (\ref{ineqms}), we get that
\begin{align*}
	\sum_{w\in\mathcal{W}}\ms((\pi_W)\pf \Si)=\ms(\Si).
\end{align*}
We have proved (\ref{ineq1}) (\ref{ineq2}) (\ref{ineq3}). Combining them gives (\ref{ineq0}) and we are done.
\end{proof}
\section{Dealing with a single stratum}\label{secss}
In this section, we will carry out the constructions outlined in Section \ref{planpf} for each separate stratum in the Almgren stratification. Our goal is to prove the following lemma, from which Theorems \ref{thmi} and \ref{thmp} follow directly in the next sections. Before stating the lemma, it is beneficial to give some definitions.% \subsection{Topological description of $X\times N_0$ and $Y\times N_f$}

The oriented submanifolds $X,Y$ are already defined in Section \ref{seceuc}. Let us recall them.

The natural splitting of $\C^{d-\dim N}$ as a Riemannian  product
\begin{align*}
	\C^{d-\dim N}=\cu{x_1\cd x_{d-\dim N}}\times \cu{y_1\cd y_{d-\dim N}},
\end{align*} 
induces a natural Riemannian product splitting on the tori
\begin{align*}
	\cz=X\times Y.
\end{align*}
From now on in this section, our ambient manifold will be
\begin{align}
	M_N=\C^{d-\dim N}/\Z^{2(d-\dim N)}\times N\times S^1=X\times Y\times N\times S^1,
\end{align}
We will reserve the labels $x,y,p,t$ for points in $X,Y,N,S^1,$ respectively. 

The projections $\pi_{X},\pi_Y,\pi_{N\times S^1}$ will denote the canonical projections of $M_N$ onto the respective factors.
\begin{defn}
	For any smooth function $$f:N\to (-\pi,\pi),$$ 
	define $N_f$ to be the oriented image of $N$ under the map $\id\times f:N\to N\times S^1,$ i.e.,
	\begin{align*}
		N_f=\{(p,f(p))|p\in N,f(p)\in(-\pi,\pi)\s S^1\}.
	\end{align*}
\end{defn}
For example, regarding $0$ as a constant function on $N,$ we get $$N_0=N\times\{0\}\s N\times S^1.$$ From now on, we will identify $N_0$ as a canonically embedded image of $N$ inside $N\times S^1.$

Now we are ready to define our integral current $\Si_N.$
\begin{defn}\label{defnsn}
	Define the integral current $\Si_N$ on $\cz\times N\times S^1$ by setting
	\begin{align*}
		\Si_N=X\times N_0+Y\times N_f.
	\end{align*}
\end{defn}
We will determine our choice of $f$ later. 

Recall Assumption \ref{assumpmff}.\begin{lem}
\label{lemr}With the same assumptions as Theorem \ref{thmi} and \ref{thmp}, for each $N\in\mff$, equip the ambient manifold
\begin{align}\label{defnmn}
M_N=\C^{d-\dim N}/\Z^{2(d-\dim N)}\times N\times S^1,
\end{align}with the product Riemannian metric, denoted by $h.$ Then for the two $d$-dimensional embedded oriented submanifolds 
\begin{align*}
X\times N_0, Y\times N_f,
\end{align*} there exists two smooth retractions 
\begin{align*}
\pi_{X\times N_0},\pi_{Y\times N_f},
\end{align*}of $M_N$ 
onto  $X\times N_0, Y\times N_f$, respectively, 
such that
\begin{enumerate}
\item In a smooth metric $g_N$ possibly different from $h,$ the smooth closed form
\begin{align}\label{eqpn1}
\phi_N=\pxn\du\dvx+\pyn\du\dvy,
\end{align}is a calibration form that calibrates the integral current $$\Si_N=X\times N_0+Y\times N_f.$$
\item The singular set of $\Si_N$ is $K_N$ in Assumption \ref{assumpmff} and is the $(\dim N)$-th stratum in the Almgren stratification of $\Si_N.$
\item For calibration mod $v$, the other three different $d$-forms
\begin{align*}	\phi_N^{+,-}=&\pxn\du\dvx-\pyn\du\dvy,\\
		\phi_N^{-,-}=&-\pxn\du\dvx-\pyn\du\dvy,\\
	\phi_N^{-,+}=&-\pxn\du\dvx+\pyn\du\dvy,
\end{align*}are also calibration forms.
\end{enumerate}%Here $X,Y$ are the subtori in $\C^{d-\dim N}/\Z^{2(d-\dim N)}$ spanned by the $x_1\cd x_{d-\dim N}$-plane and the $y_1\cd y_{d-\dim N}$, respectively. And $N_f$ is a small perturbation of coor
Here $\dvx,\dvy$ are the volume forms of $\xn,\yn$, respectively, in the $g_N$ metric.
\end{lem}
Note that the above lemma is trivially true for $\dim N=0,$ i.e., $N$ being a point, by Lemma \ref{torf}. Thus from now on we will focus on $\dim N\ge 1.$

 Roughly speaking, the proof of the above lemma consists of three steps.
\begin{enumerate}
\item In Section \ref{techlem}, we describe some basic Riemannian geometrical facts about the integral current  $$\Si_N=X\times N_0
+Y\times N_f.$$  
\item In Section \ref{secconret}, we define the retractions $\pi_{X\times N_0},\pi_{Y\times N_f}$ and the closed form \begin{align}
\phi_N=\pi_{X\times N_0}\du \dvol_{X\times N_0}^h+\pi_{Y\times N_f}\du \dvol_{Y\times N_f}^h.
\label{eqpn2}\end{align}
\item In Section \ref{secconmet}, we construct the smooth metric $g_N$ in which $\Si_N$ is calibrated by $\phi_N$ and prove that (\ref{eqpn1}) and (\ref{eqpn2}) are equivalent definitions.
\end{enumerate}
Before delving into the proof, we want to emphasize that we will encounter many different retractions and projections in this section. We adopt the following notation convention unless otherwise stated.
\begin{assump}
	The symbol $\pi$ will be reserved for the canonical projection of a product Riemannian manifold onto its factors. 
	
	The symbol $\Pi$ will be reserved for the nearest distance projection onto a Riemannian submanifold.
\end{assump}
\subsection{Basic facts about  $\Si_N=X\times N_0+Y\times N_f$}\label{techlem}
First let us describe the tangent cones to $\Si_N.$
\begin{lem}
The singular set of $\Si_N$ is the intersection of $X\times N_0$ with $Y\times N_f,$ i.e., the zero set of $f$ on $N_0$.  At any point $q\in\supp \Si_N,$  the tangent cone $C_q$ to $\Si_N$ has the following decomposition
\begin{enumerate}
\item For $q\in \reg \Si_N\cap \supp( X\times N_0),$ we have
\begin{align}\label{cpx}
C_q=\cu{x_1\cd x_{d-\dim N}}\times T_{\pi_{N\times S^1}(q)}N_0,
\end{align}
\item For $q\in \reg \Si_N\cap \supp (Y\times N_f),$ we have
\begin{align}\label{cpy}
C_q=\cu{y_1\cd y_{d-\dim N}}\times T_{\pi_{N\times S^1}(q)}N_f,
\end{align}
with $T_{\pi_{N\times S^1}(q)}N_f$ spanned by
\begin{align}\label{knf}
\ker df,\textnormal{and }\na f+|\na f|^2\pd_t,
\end{align}
and the normal space to $T_{\pi_{N\times S^1}(q)}N_f$ inside $T(N\times S^1)$ is spanned by the unit vector
\begin{align}\label{pnf}
\frac{1}{\sqrt{1+|\na f|^2}}(-\na f+\pd_t).
\end{align}
\item For $q\in \sing \Si_N$
\begin{align}\label{cps}
C_q=\cu{x_1\cd x_{d-\dim N}}\times T_{\pi_{N\times S^1}(q)}N_0+\cu{y_1\cd y_{d-\dim N}}\times T_{\pi_{N\times S^1}(q)}N_f,
\end{align}and $q$ belongs to the $(\dim N)$-th stratum in the Almgren stratification of $\Si$ when $\na f=0$ at $q$. When $\na f\not=0,$ $q$ belongs to the $(\dim N-1)$-th stratum.
\end{enumerate} 
\end{lem}
Here $\pd_t$ is the unique unit tangent vector orienting $S^1.$ %and from now on we will use coordinate label $t$ on the coordinate chart $(-\pi,\pi)$ on $S^1.$
\begin{proof}
Since $\Si_N$ is the sum of two submanifolds $X\times N_0,Y\times N_f$, let us first calculate the tangent planes to them. Due to their product structure, direct calculations give (\ref{cpx})(\ref{cpy})(\ref{knf})(\ref{pnf}).

It is clear that the singular set of $\Si_N$ must be the intersection set of $X\times N_0$ with $Y\times N_f.$ Direct calculation shows that the intersection set is the zero set of $f $ on $N_0$ identified with $N.$ Then (\ref{cps}) follows.

For a singular point $q\in \sing \Si_N,$ if $\na f=0$ at $p,$ then $T_{\pi_{N\times S^1}(q)}N_f$ coincides with $T_p N$ and thus $C_q$ has precisely $(\dim N)$-dimensional translation-invariance, which means $q$ belongs to the $(\dim N)$-th stratum in the Almgren stratification of $\Si_N.$ If $\na f\not=0,$ then $C_q$ has precisely $(\dim N-1)$-dimensional translation invariance. We are done.
\end{proof}
Secondly we need a technical lemma about exponential maps and nearest distance projections. As mentioned in Section \ref{planpf}, in general, nearest distance projections onto $X\times N_0,Y\times N_f$ may not exist everywhere. For $X\times N_0,$ this is easy to solve, as the nearest distance projection onto $X\times N_0$ coincides with the canonical projection $\pi_{X\times N}$ of $M_N$ onto the factor $X\times N$.

To deal with $Y\times N_f,$ we need the following technical lemma. Intuitively speaking, we can have uniform control on nearest distance projections of the family of submanifolds $N_{sf}$ with $s\in[-2,2]$.
\begin{lem}\label{lemexp}
Under Assumption \ref{assumpmff}, there exists a smooth function $f$ on $N,$ such that
\begin{enumerate}
\item \label{lfz}The zero set of $f$ is $f\m(0)=K_N$ and $f$ vanishes to infinite order on  $f\m(0)=K_N.$ Furthermore, $|f|<\frac{1}{2}.$
\item \label{lei}For every $s\in[-2,2],$ the closed set $N\times [-3,3]$ can be realized as the bijective image of the normal bundle exponential map $\exp^\perp_{N_{sf}}$ of $N_{sf}$ inside $N\times S^1$,
\item \label{lps}The nearest distance projection  $$\Pi_{N_{sf}}$$ of $N\times [-3,3]$ onto $N_{sf}$, regarded as a map 
\begin{align*}
[-2,2]\times N\times [-3,3]\to N\times S^1,
\end{align*}
is jointly smooth with the parameter $s$ regarded as an additional variable in $[-2,2]$.
\item On $N\times[-3,3],$ the simple form form $\Pi_{N_f}\du\dvol_{N_f}^h$ is nonzero.\label{lfn}
\end{enumerate} 
\end{lem}
Readers with experience in Riemannian geometry might find the above lemma an easy extension of the tubular neighborhood theorems. Nevertheless the author could not locate a reference, so the rest of this subsection will be devoted to the proof of the above lemma. The reader can skip this lemma on first reading through. Let us start with the construction of $f.$
\subsubsection{Construction of $f$ and definition of the extended map $E$}
It is well known that there exists a smooth function $f_{K_N}$ on $N,$ such that $f_{K_N}$ vanishes to infinite order at its zero set $f_{K_N}\m(0)=K_N.$ The central idea is to set 
\begin{align}\label{fek}
f=\e  f_{K_N}
\end{align}
 for $\e$ small enough. The smallness of $\e$ will be determined later.

Let us first discuss in detail the normal bundle exponential map of $N_{sf}$ with an additional real parameter $s\in[-2,2].$ %To facilitate calculations, it is beneficial to use the canonical parametrization of $N_{sf}$ by $N_0$.

To start, we need a way to parametrize the normal bundle of $N_{sf}$ for a fixed smooth function $f:N\to (-\pi,\pi)$ with varying $s.$ By construction, i.e., (\ref{pnf}), the normal bundle of $N_{sf}$ inside $N\times S^1$ is always trivial, i.e., diffeomorphic to $N\times \R$. 

Recall that $N_{sf}$ comes with a natural parametrization by $N$ via the map $\id\times sf.$ Using this parametrization and (\ref{pnf}), we can parametrize normal bundle of $N_{sf}$ inside $N\times S^1$  diffeomorphically (not metrically) as the product
\begin{align*}
N\times \R,
\end{align*}
via the map
\begin{align*}
	(p,\nu)\mapsto\bigg((\id\times sf)(p),\frac{\nu}{\sqrt{1+|\na (sf)|^2}}(-\na (sf)+\pd_t)\bigg).
\end{align*}
With this identification in mind, we can add $s$ to our variables parametrically as follows. 
\begin{defn}\label{defnfs}
Define a new map
$E:[-2,2]\times N\times \R\to [-2,2]\times N\times S^1$ by
\begin{align}\label{defnes}
	E(s,p,\nu)=\Bigg(s,\exp\bigg((\id\times sf)(p),\frac{\nu}{\sqrt{1+|\na (sf)|^2}}(-\na (sf)+\pd_t)\bigg)\Bigg).
\end{align} %\begin{align}E(s,p,\nu)=\bigg(s,\Big(\exp^\perp_{N_{sf}}\circ \big(\id\times sf\big)\Big)(p,\nu)\bigg).\end{align}
\end{defn}
%Here $\exp^\perp_{N_{sf}}$ is the normal bundle exponential map of $N_{sf}.$
Here $$\exp(q,v)$$ denotes image of the exponential map of a point $(q,v)$ in the tangent bundle $$T(N\times S^1)=\{(q,v)|q\in N\times S^1,v\in T_q(N\times S^1)\}$$ of $N\times S^1.$ %Then using , we can rewrite (\ref{defnes}) as
\subsubsection{Estimating $dE$}
From the definition of Jacobi fields as variation fields of geodesics (\cite[Chapter 10]{JLr}), it is straightforward to see that
at a point $(q,v)$, the differential of the exponential map at a tangent vector $(\de q,\de v)$ to the double tangent bundle $T_{(p,v)}T(N\times S^1)$ is \begin{align*}
d\exp(\de q,\de v)=J(1),
\end{align*}where $J(1)$ is the Jacobi field at time $1$ along the geodesic $\exp(q,\tau v)$ for $\tau\in[0,1]$, with initial values $J(0)=\de q,J'(0)=\de v.$

From now on, we set $J(1)$ to be the Jacobi field at time $1$ along the geodesic $\exp(p,\tau \nu)$ for $\tau\in[0,1]$, with 
\begin{align*}
	J(0)=f(p)\pd_t,\\J'(0)=-\nu\na f.
\end{align*}Integrating the Jacobi equation twice and using Gronwall's inequality, we get that 
\begin{align}\label{gronw}
	\no{J(1)}=O((1+|\nu|)\no{f}_{C^1}),
\end{align}with big $O$ depending only on the ambient metric $h.$
 
On the other hand, direct calculation using chain rule shows that at any point $$(0,p,\nu)\in\{0\}\times N\times \R,$$ we have %1at $T_{(0,\exp(p,\nu))}([-2,2]\times N\times S^1),$
\begin{align}
\label{des}dE(\pd_s)=&(\pd_s,J(1)),\\\label{dep}dE(\de p)=&(0,\de p),\\\label{dev}dE(\pd_\nu)=&(0,\pd_t).
\end{align}
Here $\de p $ is any vector tangent to $N$ in the natural  Riemannian splitting of $T(N\times S^1)\cong TN\times TS^1.$ 

Thus, by setting $\e$ small in (\ref{fek}) and restrict to $\nu\in[-3,3]$, using (\ref{gronw}) we can make sure that $$\no{J(1)}<\frac{1}{100}.$$ Putting this into (\ref{des}) and considering (\ref{dep}) (\ref{dev}), we deduce that $dE$ is a vector space isomorphism at any point $(0,p,\nu)\in\{0\}\times N\times \R,$ with uniform control on the representing matrix of $dE$ in coordinate.

Now using the quantitative version of implicit function theorem (\cite[Section 8]{MC}), we deduce that 
\begin{fact}\label{fctinv}
There is a radius $r>0$ such that for any $$(0,p,\nu)\in \{0\}\times N\times [-3.14,3.14],$$ the map $E$ restricted to the radius $r$ geodesic ball centered at $(0,p,\nu)$ on $[-2,2]\times N\times \R$ is a diffeomorphism.
\end{fact}
\subsubsection{The map $E$ is a diffeomorphism for $s$ small}
With Fact \ref{fctinv} in hand, we will prove that
\begin{fact}\label{fctdiff}
There is $s_0>0,$ such that $E$ restricted to $$[-s_0,s_0]\times N\times[-3.14,3.14],$$
is a diffeomorphism.
\end{fact}
By \cite[Proposition 4.22(c)]{JL}, to prove Fact \ref{fctdiff}, it suffices to show that there is $s_0>0$ such that $E$ restricted to $[-s_0,s_0]\times N\times[-3.14,3.14]$ is an injective immersion. By Fact \ref{fctinv}, $E$ is an immersion in  $[-s_0,s_0]\times N\times[-3.14,3.14]$ for $s_0>0.$ It suffices to verify the injectivity of $E$ in $[-s_0,s_0]\times N\times[-3.14,3.14]$ with a possibly smaller $s_0>0.$

Suppose the contrary. In other words, suppose there exists a sequence of positive real numbers $\{s_j\}$ with $\lim_{j\to\infty}s_j=0,$ and two sequence of points $$\{(s_j,p_j,v_j)\},\{(s_j,p_j',v_j')\}\s [-2,2]\times N\times [-3.14,3.14],$$
such that 
\begin{align*}
E(s_j,p_j,v_j)=E(s_j,p_j',v_j'),
\end{align*}for all $j.$
By compactness of $[-2,2]\times N\times [-3.14,3.14]$, we can pass to two subsequences, not relabeled, such that
\begin{align*}
\lim_{j\to\infty}p_j=p,\lim_{j\to\infty}p_j'=p',
\lim_{j\to\infty}v_j=v,\lim_{j\to\infty}p_j'=v'.
\end{align*}
This implies that
\begin{align}\label{eqe}
E(0,p',v')=E(0,p,v).
\end{align}
The restriction of $E$ to $\{s=0\}$ constant is exactly the normal bundle exponential map of $N_{0}$. It is straightforward to check that the normal bundle exponential map of $N_0$ is injective for all normal vectors with length less than $\pi.$  Thus, (\ref{eqe}) implies that $p'=v',p=v.$

Consequently, for $j$ large enough, $p_j,p_j',v_j,v_j'$ stays in a geodesic ball of radius $r$ centered at some point on $\{0\}\times N\times[-3.14,3.14].$ However, this is a contradiction to Fact \ref{fctinv}. Thus, we have finished the proof of Fact \ref{fctdiff}.
\subsubsection{Wrapping up the proof of Lemma \ref{lemexp}}
Now apply a change of variable $s\to \frac{s_0}{2}s$, then Fact \ref{fctdiff} holds with $2$ instead of $s_0.$ This corresponds to setting a new $f\to \frac{s_0}{2}f.$ 

The restriction of $E$ to each $\{s\}\times N\times[-3.14,3.14]$ is exactly the normal bundle exponential map of $N_{sf}$ for vectors with length at most $3.14.$ The product structure of $N\times S^1$ implies that projection of geodesics onto $S^1$ will be constant speed curves. Using (\ref{pnf}), we deduce that $N\times[-3,3]$ lies in the diffeomorphic image of $\{s\}\times N\times[-3.14,3.14]$ for all $s\in[-2,2],$ provided $\no{f}_{C^1}$ is small. This smallness of $\no{f}_{C^1}$ is obtained by setting $\e$ in (\ref{fek}) small.

To wrap it up, we have proved the second bullet (\ref{lei}) of our lemma. The first bullet (\ref{lfz}) was already proven when choosing $f.$

To prove the third bullet (\ref{lps}), note that whenever the normal bundle exponential map $\exp^\perp_{N_{sf}}$ is bijective, we can identify the nearest distance projection onto $N_{sf}$ as the projection of the normal bundle of $N_{sf}$ onto its base, i.e., $N_{sf}$ itself. %by (\ref{fctdiff}), the normal bundle exponential map of $N_{sf}$ on normals with length at most $3.14,$ i.e., the restriction of $E$ to $s\times N\times[-3.14,3.14]$, is a diffeomorphism. However,
In other words, we have
\begin{align*}
\Pi_{N_{sf}}=(\id\times sf)\circ\pi_N\circ(E\m)\circ i_s,
\end{align*}where $i_s:N\times S^1\to [-2,2]\times N\times S^1$ is the map $i_s(p,t)=(s,p,t)$, and $\pi_N$ is the canonical projection of $[-2,2]\times N\times \R$ onto $N.$ Since every term in the composition is jointly smooth with $s$ as an additional variable. We deduce that $\Pi_{N_{sf}}$ is jointly smooth with $s$ added as a variable.

To prove the fourth bullet (\ref{lfn}), %note that $N\times[-3,3]$ is contained in the bijective image of the normal bundle exponential map of $N_f$ on normals with length at most $3.14.$ Thus, by  identifying
by bullet (\ref{lei}), we can identify $N\times[-3,3]$ with a closed set with non-empty interior in the normal bundle of $N_f$ restricted to normals of length at most $3.14,$ i.e., a closed set with non-empty interior in the trivial bundle $$N_f\times [-3.14,3.14].$$ However, the nearest distance projection under this identification is precisely the canonical projection of $N_f\times [-3.14,3.14]$  onto the $N_f$ factor. Thus, $\Pi_{N_f}\du \dvol^h_{N_f}$ in $N_f\times[-3.14,3.14]$ is simply the constant extension of the non-zero form $\dvol^h_{N_f}$ along the $[-3.14,3.14]$ directions, thus nonzero.
\subsection{Construction of the retractions}\label{secconret}
From now on, we always use $f$ as constructed in Lemma \ref{lemexp}.
We now define a retraction $\pi_{N_f}:N\times S^1\to N_f$.
\begin{defn}\label{defnpnf}Let $\ka:\R\to\R$ be a smooth monotonically decreasing function that equals to $1$ on $(-\infty,2.25]$ and equals to $0$ on $[2.75,\infty).$
Define
\begin{align*}
\pi_{N_f}(p,t)=\begin{cases}
\Pi_{N_f}(p,t),&\textnormal{for }|t|< 2,\\
(\id\times f)\circ\pi_N\circ\Pi_{N_{\ka(|t|)f}}(p,t),&\textnormal{for }2\le|t|\le 3,\\
(\id\times f)\circ\pi_{N}(p,t),&\textnormal{for }|t|>|3|,
\end{cases}
\end{align*}
\end{defn}Here we regard $S^1$
as $[-\pi,\pi]$ with two end points identified and $\pi_N$ is the canonical projection of $N\times S^1$ onto $N$. Using Lemma \ref{lemexp} and the fact that 
\begin{align*}
	(\id\times f)\circ\pi_N\circ \Pi_{N_f}=\Pi_{N_f},\\
	\pi_N\circ\Pi_{N_0}=\pi_N^2=\pi_N,
\end{align*} it is straightforward to verify that $\pi_{N_f}$ is a smooth retraction of $N\times S^1$ onto $N_f.$

Intuitively speaking, $\pi_{N_f}$ transits from $\Pi_{N_f}$, which is not defined everywhere, to $(\id\times f)\circ\pi_N$, which is defined everywhere, as $|t|$ goes from $2$ to $3$.

Use $(x,y,p,t)$ to denote points in 
\begin{align*}
M_N=X\times Y\times N\times S^1,
\end{align*}
i.e., $x,y,p,t,$ respectively are the canonical projected images onto the factors $X,Y,N,S^1$ of $M_N,$ respectively. 

Now we are ready to define our retractions  $M_N$ onto $X\times N_0,Y\times N_f.$
\begin{defn}\label{defnp}
Define twos maps $\pxn,\pyn:M_N\to X\times N_0$ by
\begin{align*}
\pxn(x,y,p,t)=&(x,0,p,0),\\
\pyn(x,y,p,t)=&(0,y,\pi_{N_f}(p,t)),
\end{align*}
or in notation of projections
\begin{align*}
\pxn=&(\pi_X,0,\pi_N,0),\\
\pyn=&(0,\pi_Y,\pi_{N_f}\circ\pi_{N\times S^1}),
\end{align*}
where $\pi_X,\pi_Y,\pi_{N\times S^1},$ respectively are the canonical projection of $M_N$ on to the factors $X,Y,N\times S^1$, respectively.
\end{defn}
By construction, $\pxn,\pyn$ are smooth retractions onto $X\times N_0,Y\times N_f$, respectively.
\subsection{The calibration form and the metric}\label{secconmet}
\begin{defn}
Define four $d$-forms
\begin{align*}
\phi_N=&\pxn\du\dvol^h_{X\times N_0}+\pyn\du\dvol^h_{Y\times N_f},\\	\phi_N^{+,-}=&\pxn\du\dvol^h_{X\times N_0}-\pyn\du\dvol^h_{Y\times N_f},\\
	\phi_N^{-,+}=&-\pxn\du\dvol^h_{X\times N_0}+\pyn\du\dvol^h_{Y\times N_f},\\
	\phi_N^{-,-}=&-\pxn\du\dvol^h_{X\times N_0}-\pyn\du\dvol^h_{Y\times N_f}.
\end{align*}
Set
\begin{align}
\label{mnf}\ms(\pi_{N_f})=\max\bigg\{&1,(\cms_h\phi_N)^{\frac{2}{d}},(\cms_h\phi^{+,-}_N)^{\frac{2}{d}},(\cms_h\phi^{-,+}_N)^{\frac{2}{d}},(\cms_h\phi^{-,-}_N)^{\frac{2}{d}},\\&\sup_{X\times Y\times N\times[-2,2]}\no{\pi_{N\times S^1}\du\pi_{N_f}\du \dvol_{N_f}^h}_h^{-\frac{2}{d-\dim N}},\\&\sup_{X\times Y\times N\times[-2,2]}\no{\pi_{N\times S^1}\du\pi_{N_f}\du \dvol_{N_f}^h}_h^{\frac{2}{\dim N}}\bigg\}.\end{align}
\end{defn}
%Here $\cms_{h,N\times[-2,2]}\pi_{N_f}\du \dvol_{N_f}^h$ means the comass of $\du \dvol_{N_f}$ in the Riemannian  manifold $N\times[-2,2]\s N\times S^1.$
Note that the above definition of forms is slightly different from the one in Lemma \ref{lemr}. This is not a mistake, as we will show later that $\dvol^h_{X\times N_0}=\dvx$ and $\dvol^h_{Y\times N_f}=\dvy.$

 The Riemannian product structure of 
\begin{align*}
M_N=X\times Y\times(N\times S^1),
\end{align*}
provides a natural orthogonal splitting of the base metric $h$ into components
\begin{align*}
h=h_X+h_Y+h_{N\times S^1}.
\end{align*}
\begin{defn}\label{defngn}
Define an auxiliary metric $h'$ on $M_N$ by
\begin{align*}
h'=\no{\pi_{N\times S^1}\du\pi_{N_f}\du \dvol_{N_f}^h}_h^{-\frac{2}{d-\dim N}}h_X+h_Y+\no{\pi_{N\times S^1}\du\pi_{N_f}\du \dvol_{N_f}^h}_h^{\frac{2}{\dim N}}h_{N\times S^1}.
\end{align*}
Define a metric $g_N$ on $M_N$ by
\begin{align*}
g_N=\begin{cases}
h',&\textnormal{for }|t|\le 1,\\
\ka(|t|-1)h'+\big(1-\ka(|t|-1)\big)\ms(\pi_{N_f})h,&\textnormal{for }1<|t|\le 2,\\
\ms(\pi_{N_f})h,&\textnormal{for }|t|> 3,
\end{cases}
\end{align*}
\end{defn}
Here $\ka$ is defined in Definition \ref{defnpnf} and $\pi_{N\times S^1}$ is the canonical projection of $M_N$ onto $N\times S^1.$ 

The rest of this section will be devoted to the proof of the following lemma, which immediately gives us Lemma \ref{lemr}.
\begin{lem}\label{lemfff}We have the following.
\begin{enumerate}
\item The metric $g_N$ is smooth and the forms $\phi_N,\phi_N^{+,-},\phi_N^{-,+},\phi_N^{-,-}$ are closed and smooth. 
\item Both  $\pxn\du\dvol^h_{X\times N_0}$ and $\pyn\du\dvol^h_{Y\times N_f}$ are unit simple forms in $h'$ for $|t|\le 2$ and they satisfy
\begin{align*}
	\pxn\du\dvol^h_{X\times N_0}=&\pxn\du\dvx,\\
	\pyn\du\dvol^h_{Y\times N_f}=&\pyn\du\dvy.
\end{align*} \label{cmsh'}
\item $\phi_N,\phi_N^{+,-},\phi_N^{-,+},\phi_N^{-,-}$ are calibration forms in metric $h'.$\label{ffmh'}
\item $\phi_N,\phi_N^{+,-},\phi_N^{-,+},\phi_N^{-,-}$ are calibration forms in metric $g_N.$
\item $\phi_N$ calibrates $\Si_N$ in $g_N.$
\end{enumerate}

\end{lem}
%In this, we will verify the first bullet and the next subsections will be devoted to the proof the second and the third bullet.
%Let us start with the smoothness of $\phi_N$ and $g_N.$
\subsubsection{Smoothness of $\phi_N$ and $g_N$}
By construction, $\pi_{X\times N_0},\pi_{Y\times N_f}$ are smooth, from which the smoothness of $\phi_N^{+,-},\phi_N^{-,+},\phi_N^{-,-}$ follows. %It is also evident that $\phi_N^{+,-},\phi_N^{-,+},\phi_N^{-,-}$ are closed forms, as they are pullbacks of volume forms of smooth submanifolds, which are closed by definition.

For the smoothness of $g_N,$ let us first note that the Riemannian length $$\no{\pi_{N_f}\du \dvol^h_{N_f}}_h$$ has positive lower bound on $X\times Y\times N\times [-2,2]$ by bullet (\ref{lfn}) in Lemma \ref{lemexp} and Definition \ref{defnpnf}. 

Thus, we deduce that
both $\no{\pi_{N_f}\du \dvol^h_{N_f}}_h^{\frac{2}{\dim N}}
$ and $\no{\pi_{N_f}\du \dvol^h_{N_f}}_h^{-\frac{2}{d-\dim N}}
$ are smooth, bounded and positive on $N\times[-2,2].$  The Riemannian product structure of $M_N=X\times Y\times (N\times S^1)$ implies that
$$\no{\pi_{N\times S^1}\du\pi_{N_f}\du \dvol_{N_f}^h}_h=\no{\pi_{N_f}\du \dvol_{N_f}^h}_h\circ\pi_{N\times S^1}.$$This implies that $\ms(\pi_{N_f})$ in (\ref{mnf}) is well-defined and finite and the auxiliary metric $h'$ is smooth and well-defined on $N\times[-2,2].$ The smoothness of $g_N$ follows immediately.

The closedness of the four differential forms follows from the closedness of volume forms and the fact that pullback preserves closedness.
\subsubsection{Comass of $\phi_N$ in $h'$}
First of all, by (\ref{lfz}) of Lemma \ref{lemexp}, and Definition \ref{defnsn},  $\Si_N$ is supported on $X\times Y\times N\times[-\frac{1}{2},\frac{1}{2}]$.

Now let us prove bullet (\ref{cmsh'}) of Lemma \ref{lemfff}.

The product structure of $M_N$ gives a natural orthogonal splitting 
\begin{align}
	\label{tms}TM_N=&TX\oplus TY\oplus T(N\times S^1),\\
		\label{tmds}T\du M_N=&T\du X\oplus T\du Y\oplus T\du (N\times S^1).
\end{align} Note that
\begin{fact}
Both	 $h'$ and $g_N$ both respect the product structure of $M_N$, so the above orthogonal splitting holds in $h,h'$ and $g_N.$
\end{fact}

 On $X\times Y\times N\times[-2,2],$ using the Riemannian product structure of $M_N$, we deduce from Definition \ref{defnp} that
\begin{align}\label{pxn0}
	\pxn\du\dvol_{X\times N_0}^h=\pi_X\du\dvol_X^h\w \pi_N\du\dvol_{N}^h.
\end{align}
Apply (\ref{eqw}) of Lemma \ref{lemvfm}, we deduce that
\begin{align}\label{cmsx}
\cms_{h',X\times Y\times N\times[-2,2]}	\pxn\du\dvol_{X\times N_0}^h=1.
\end{align}
 Straightforward calculation using (\ref{eqn})  of Lemma \ref{lemvfm}, with (\ref{cpx}) and (\ref{cpy}) in mind, shows that \begin{align}
\label{pxn01}	\pxn\du\dvol_{X\times N_0}(T^{h'}_{(x,0,p,0)}(X\times N_0))=1,\\
\label{pxn00}	\pxn\du\dvol_{X\times N_0}(T^{h'}_{(0,y,p,t)}(Y\times N_f))=0,
\end{align} 
where $T^{h'}_{(x,0,p,0)}(X\times N_0)$ is a unit $d$-vector orienting the tangent space to $X\times N_0$ and $T^{h'}_{(0,y,p,t)}(Y\times N_f)$ is a unit $d$-vector orienting the tangent space to $Y\times N_f.$ Note that (\ref{pxn01}) also implies that $$\dvol^h_{X\times N_0}=\dvx,$$ since $h'=g_N$ on $X\times N_0$ by Definition \ref{defngn} and $\pxn$ is a retraction onto $X\times N_0.$

Similar calculation shows that in metric $h'$ on $X\times Y\times N\times[-2,2],$ we have the decomposition
\begin{align*}
	\pyn\du\dvol_{Y\times N_f}^h=\pi_Y\du\dvol_{Y}^h\w\pi\du_{N\times S^1}\pi_{N_f}\du\dvol_{N_f}^h.
\end{align*}
Apply (\ref{eqw}) of Lemma \ref{lemvfm}, we deduce that
\begin{align}\label{cmsy}
	\cms_{h',X\times Y\times N\times[-2,2]}	\pyn\du\dvol_{Y\times N_f}^h=1.
\end{align}
Straightforward calculation using (\ref{eqn})  of Lemma \ref{lemvfm}, with (\ref{cpx}) and (\ref{cpy}) in mind, shows that 
\begin{align}
	\pyn\du\dvol^h_{Y\times N_f}(T^{h'}_{(x,0,p,0)}(X\times N_0))=0,\\
\label{pyn01}	\pyn\du\dvol^h_{Y\times N_f}(T^{h'}_{(0,y,p,t)}(Y\times N_f))=1.
\end{align} 
This time (\ref{pyn01}) implies that $$\dvol^h_{Y\times N_f}=\dvy$$ as $h'=g_N$ on $X\times N_0$ by Definition \ref{defngn}, and $\pyn$ is a retraction onto $Y\times N_f.$

Using (\ref{cmsx}) and (\ref{cmsy}), by Lemma \ref{torf}, we deduce (\ref{ffmh'}) of Lemma \ref{lemfff}.
\subsubsection{Computing the comass of $\phi_N$ in $g_N$}
On $X\times Y\times N\times[-1,1]$, we have $g_N=h'.$ Since $\Si_N$ is supported in $X\times Y\times N\times[-1,1]$, we deduce that $\phi_N$ calibrates $\Si_N$ in $X\times Y\times N\times[-1,1]$ by (\ref{pxn01}) and (\ref{pyn01}).

Using (\ref{cms1}) of Fact \ref{cmsvec} and $g_N\ge h'$ as quadratic forms on $X\times Y\times N\times [-2,2].$ we verify that $\phi_N$ has comass at most $1$ on $X\times Y\times N\times [-2,2]$ by (\ref{ffmh'}) of Lemma \ref{lemfff}.

Using (\ref{cms2}) of Fact \ref{cmsvec}, we deduce that $\phi_N$ has comass at most $1$ on $X\times Y\times N\times \{|t|\ge 2\}.$

The same argument works for the other three different forms $\phi_N^{+,-},\phi_N^{-,+},\phi_N^{-,-}.$ This finishes the proof of Lemma \ref{lemfff}, from which Lemma \ref{lemr} follows.
\section{Proof of Theorem \ref{thmi}, and Theorem \ref{thms}}\label{secpfmt}
In this section, we will prove all the theorems.
\subsection{Proof of Theorem \ref{thmi}}\label{secpthmi}
Recall Lemma \ref{lemr} and Assumption \ref{assumpmff}. For each $N\in \mff$, set $$M_N'=M_N\times (S^1)^{\dim N-k},$$ where we regard $(S^1)^{0}$ as a single point. Apply Lemma \ref{lemprod} to $M_N'$. Then apply Lemma \ref{lemdis} to the disjoint union $$\bigcup_{N\in\mff}M_N'.$$ We are done.
\subsection{Proof of Theorem \ref{thmp}}
Apply Lemma \ref{lemcalv} to Lemma \ref{lemr}. Then argue as in the last subsection using Lemma \ref{lemdis} and Lemma \ref{lemprod}.
\subsection{Proof of Theorem \ref{thms}}
The proof is very similar to that of Theorem \ref{thmi}. We only give a sketch.  Then the only part of Theorem \ref{thmi} that is not covered yet is to deal with the $(d-1)$-stratum. In this case, Lemma \ref{torf} is not applicable, as both $X$ and $Y$ are of dimension $1$. However, the analogues of Lemma \ref{lemr} hold, if instead of requiring $\phi_N=\pxn\du\dvx+\pyn\du\dvy$ to be a calibration form, we require the two factors $\pxn\du\dvx$ and $\pyn\du\dvy$ both to be calibration forms separately.  The last bullet folows because both $X\times N_0$ and $Y\times N_f$ are homologically area-minimizing, and thus for diffeomorphisms $\Phi$ homotopic to the identity we have, \begin{align*}
	\ms(\Phi\pf(X\times N_0+Y\times N_f))=&\ms(\Phi\pf (X\times N_0))+\ms(\Phi\pf(Y\times N_f))\\\ge& \ms(X\times N_0)+\ms(Y\times N_f)\\=&\ms(X\times N_0+Y\times N_f).
\end{align*}  Consequently, the chain sum of $X\times N_0$ and $Y\times N_f$ is a stable stationary varifold. Now use Lemma \ref{lemdis} and \ref{lemprod} on $X\times N_0$ and $Y\times N_f$ separately for all $N\in\mathcal{F}$. Argue as in the previous subsections. We are done.
\section{Remarks and discussions}\label{conc}
\subsection{Comparison with Professor Leon Simon's work}\label{ls}
Professor Leon Simon (\cite{LSfr}) has done the groundbreaking construction of stable stationary hypervarifolds of dimension $d\ge 8$ and codimension $c=1,$ with the singular set being any closed subset of $\R^{d-7}$. Indeed, Professor Simon's work is one of the main motivations for the author to pursue this problem. Just for the purpose of the peer reviewing process, we comment on the differences between the two. First, it is not known if his examples are area-minimizing.  Second, the results in his work are based on codimension $1$ area-minimizing cylinders, while ours are of codimension $c\ge 3$, thus permitting relatively larger singular sets. Third, the singular set of our area-minimizing submanifolds comes from self-intersection, or more intuitively, of crossing type, while his singular set comes from smoothing of conical points, thus in some sense a more genuine singular set. Finally, our Theorem \ref{thms} can be seen as an analogue of his results in higher codimension.
\subsection{Why the metrics are so complicated in Section \ref{secss}}
In general the comass of a form is not equal to the Riemannian length of the form. Furthermore, the \textbf{pointwise} comass of the form depends only continuously but not necessarily smoothly on the form and the metric. Thus, many prevalent constructions, like normalizing, do not work at all for comass. For instance, normalizing a non-vanishing form by its  \textbf{pointwise} Riemannian length always yields a smooth form of unit length, but normalizing a non-vanishing form by its \textbf{pointwise} comass can yield a non-smooth form. The reader can take (\ref{cmsc}) for example, where the Lipschitz operator $\max$ appears naturally in the value of comass. For another example, classifying non-vanishing constant forms of unit length on Euclidean space is a trivial task. However, tremendous efforts in the 1980s only {partially} classified constant forms of unit comass in \cite{DHM} up to forms in $\R^8,$ and there is no further progress in the {past 40 years}! This illustrates the difficulty of working with comass, and why we often have to split the ambient metrics into different components and rescale differently.

Unfortunately, the comass of differential forms is canonically the dual norm (\cite{HFrf}) to the mass of currents. Thus, there is no way to replace comass with other simpler norms \textbf{in general}. This might partially explain why no previous work has succeeded in finding area-minimizing currents with fractal singular sets.
\subsection{Type of singularities we use}
	The singularities we use all come from non-transverse self-intersections. It is unclear at first sight whether under generic perturbations of ambient metrics, the fractal singular sets will change. We will provide a partial solution in \cite{ZLa2}. The generic transversality results of Brian White \cite{BWgt,MR1305283} does not apply directly, since there is no guarantee that area-minimizers in nearby metrics are minimal immersions. However, by construction, there are non-generic ways to dissolve the singularities by perturbing $f$ in Section \ref{secss}.
\subsection{The case of analytic metrics}
	In real analytic metrics, immersed minimal submanifolds are locally real analytic subvarieties, so they cannot possess fractal self-intersections. Thus, if we strengthen Almgren's conjecture \cite[ Problem 5.4]{GMT} by imposing the restriction of real analytic metrics, then the stronger statement still remains open, no matter how we adapt our method. In particular, we still do not know even if stationary minimal surfaces can possess a fractal singular set in Euclidean space with the standard flat metric. Professor Simon's work \cite{LSfr} also needs to perturb the metric  non-analytically. Moreover, around any singular point, our construction is locally reducible into two smooth pieces. It would be very interesting to get locally non-irreducible examples, for example, with true branch points being the singular set. In some sense, the locally non-irreducible singularities are the more genuinely singular points. One can compare the situation in \cite[Theorem 1.8 and Conjecture 1.7 ]{DPHM}.

Finally, we want to point out a heuristic towards the existence of exotic singular sets in analytic metrics. Thom's classical results \cite{RT} show that general homology classes cannot be represented by maps from smooth manifolds. On the other hand, by the resolution of singularity for subanalytic subvarieties \cite{BM}, this implies that general homology classes cannot be represented by subanalytic varieties. Since in analytic metrics, the regular sets of area-minimizing integral currents are always locally analytic subvarieties, it would be very weird if the singular set per se is subanalytic for area-minimizing integral currents in these general homology classes with no subanalytic representatives. However, the reader should keep \cite{SA} in mind as a counterargument to the above heuristic.
\subsection{Different strata interacting with each other}
Yang Li and Yongsheng Zhang have raised the question of whether different strata in the Almgren stratification of area-minimizing current with fractal singular sets can limit to each other, instead of being a disjoint union as presented here. The answer is yes. We will pursue it in \cite{ZLa2}.
\subsection{Whether this happens on general manifolds}
The constructions in this manuscript can be modified to produce $d$-dimensional area-minimizing integral currents with fractal singular sets on any manifold $M^{d+c}$ with $H_d(M,\R)\not=0$ and $d,c\ge 3$. We will pursue this  in \cite{ZLa2}. Thus, area-minimizing integral currents with fractal singular sets do appear on general manifolds.

On the other hand, it is easy to see that if we set $N=(S^1)^n$ in Section \ref{secss}, and let $f$ have an arbitrarily small $C^\infty$ norm, we can obtain area-minimizing integral currents with fractal singular sets in metrics arbitrarily close to flat metrics on tori. Thus, area-minimizing integral currents with fractal singular sets appear naturally even in metrics arbitrarily close to canonical space form metrics.
\subsection{History of constructions}	One might ask why our method cannot give singularities in the top-dimensional stratum, i.e., branch/cusp-like points. The idea is that the tangent cones at the singular points of the currents we use, lie in the interior of the moduli space of pairs of area-minimizing planes. The branch points correspond to the degenerate boundary points of the moduli space, thus having fewer degrees of freedom to perturb. The calibrations in this paper are inspired by the work of Frank Morgan (\cite{FMtd},\cite{FMeu}), Gary Lawlor (\cite{GL}) and Dana Nance (\cite{DN}) on the angle conjecture and intersecting minimal surfaces. 
	\printbibliography
\end{document}